\documentclass[12pt,a4paper,leqno,titlepage]{amsart}

\usepackage{amsmath,amscd,amsfonts,amsthm,amssymb,dsfont,bbm,manfnt,amsaddr}
\usepackage{hyperref,url,color,appendix,eucal}
\usepackage[usenames,dvipsnames]{xcolor}
\usepackage{multicol}
\usepackage[english]{babel}
\usepackage[OT2,T1]{fontenc}
\usepackage[utf8]{inputenc}
\usepackage[all]{xy}
\usepackage{graphicx}

\hypersetup{colorlinks=true,linkcolor=black,citecolor=black,urlcolor=Periwinkle,pdfborderstyle={/S/U/W 1}}
\newenvironment{pf}
{\medskip\noindent {\it Proof.  }}
{\hfill\nobreak $\Box$ \par\bigbreak}

\newcommand{\isomo}{\overset{\sim}{\rightarrow}}

\newcommand{\GL}{\mathrm{GL}}

\newcommand{\ps}{\par \smallskip}

\newcommand{\Z}{\mathbb{Z}}

\newcommand{\Q}{\mathbb{Q}}

\newcommand{\R}{\mathbb{R}}    
\newcommand{\C}{\mathbb{C}}

\newtheorem{thm}[subsection]{Theorem}
\newtheorem*{thm*}{Theorem}
\newtheorem{lemme}[subsection]{Lemma}

\newtheorem*{lemme*}{Lemma}
\newtheorem{remark}[subsection]{Remark}

\newtheorem{cor}[subsection]{Corollary}

\newtheorem{prop}[subsection]{Proposition}

\newtheorem{example}[subsection]{Example}

\newtheorem{definition}[subsection]{Definition}

\newtheorem{thmintro}{Theorem}

\makeatletter

\@addtoreset{equation}{section}
\makeatother

\hypersetup{colorlinks=true,linkcolor=black,citecolor=black,urlcolor=Periwinkle,pdfborderstyle={/S/U/W 1}}

\title[]{{\sc An automorphic generalization of the Hermite-Minkowski theorem}}

\author{Ga\"etan Chenevier}

\thanks{\noindent  L'auteur a \'et\'e financ\'e par le C.N.R.S. et a re\c{c}u le soutien du projet ANR-14-CE25 (PerCoLaTor).}

\address{{\small  C.N.R.S., Universit\'e Paris-Sud, Institut de Math\'ematiques d'Orsay, \\ B\^atiment 307, 91405 Orsay, France,\\ 
{\rm \href{mailto:gaetan.chenevier@math.cnrs.fr}{gaetan.chenevier@math.cnrs.fr}}}}

\date\today
\begin{document}

\begin{abstract} We show that for any integer $N\geq 1$, there are only finitely many cuspidal algebraic automorphic representations of $\GL_n$ over $\Q$, with $n$ varying, whose conductor is $N$ and whose weights are in the interval $\{0,\dots,23\}$. \ps
More generally, we define an explicit sequence $({\rm r}(w))_{w\geq 0}$ such that for any number field $E$ whose root-discriminant is $<{\rm r}(w)$, and any ideal $\mathcal{N}$ in the ring of integers of $E$, there are only finitely many cuspidal algebraic automorphic representations of $\GL_n$ over $E$, with $n$ varying, whose conductor is $\mathcal{N}$ and whose weights are in the interval $\{0,\dots,w\}$. We also show that, assuming a version of {\rm GRH}, we may replace ${\rm r}(w)$ with $8 \pi e^{-\psi(1+w)}$ in this statement. \ps
The proofs are based on some new positivity properties of certain real quadratic forms which occur in the study of the Weil explicit formula for Rankin-Selberg ${\rm L}$-functions. Both the effectiveness and the optimality of the methods are discussed.  \end{abstract}

{\let\newpage\relax\maketitle}

\section{Introduction}

A classical result in the geometry of numbers, due to Hermite and Minkowski, asserts that there are only finitely many number fields with given discriminant. Minkowski's contribution, proved using his geometry of numbers, is that such number fields have a bounded degree. These results extend to Artin representations. Indeed, as is shown in \cite{abcz} there are only finitely many isomorphism classes of complex representations of ${\rm Gal}(\overline{\Q}/\Q)$ with given conductor and dimension; according to Odlyzko \cite{odlyzkocond}, this even holds without fixing the dimension if we assume the Artin conjecture. In a similar spirit, famous results of Faltings and Zahrin show that there are only finitely many isomorphism classes of abelian varieties over $\Q$ with given conductor and dimension; Mestre proved that this still holds without fixing the dimension, assuming standard conjectures about the ${\rm L}$-functions of the Tate modules of abelian varieties over $\Q$ \cite[\S III]{mestre}. \ps

In this paper, we prove a new finiteness result for algebraic automorphic representations of ${\rm GL}_n$ over number fields, which is reminiscent of these contributions of Minkowski, Odlyzko and Mestre. In very special cases, it actually explains the finiteness statements above. \ps


\subsection{Statements} Let $E$ be a number field and $\pi$ a cuspidal automorphic representation of $\GL_n$ over $E$. Following Harish-Chandra and Langlands, the infinitesimal characters of the Archimedean components of $\pi$ may be viewed as a collection of semisimple conjugacy classes in ${ \rm M}_n(\C)$ (\S \ref{archalg}); the eigenvalues of those classes will be called the {\it weights} of $\pi$. The representation $\pi$ is said {\it algebraic} if its weights are in $\Z$; its {\it motivic weight} ${\rm w}(\pi)$ is then defined as the sum of the smallest and of the largest weights of $\pi$. As twisting a cuspidal $\pi$ by $|\det |^{\pm 1}$ simply shifts the weights by $\pm 1$, we will often assume that the weights of an algebraic $\pi$ are nonnegative. In dimension $n=1$, an algebraic $\pi$ is a Hecke character {\it of type} ${\rm A}_0$ in the sense of Weil. In dimension $2$ and say with $E=\Q$, an algebraic $\pi$ with weights $0$ and ${\rm w}(\pi) \geq 0$ is generated either by a cuspidal elliptic eigenform of {\it usual weight} ${\rm w}(\pi)+1$, or by a Maass form with Laplacian eigenvalue $\frac{1}{4}$ in the case ${\rm w}(\pi)=0$. \ps\ps

By classical conjectures of Fontaine-Mazur \cite{fm} and Langlands \cite{langlandsedin}, the cuspidal algebraic automorphic representations of ${\rm GL}_n$ are of considerable interest in arithmetic geometry: they are expected to correspond bijectively to isomorphism classes of {\it geometric} irreducible representations $\rho : {\rm Gal}(\overline{E}/E) \rightarrow {\rm GL}_n(\overline{\Q_\ell})$ (or motives). In this correspondence, the Godement-Jacquet ${\rm L}$-function of a $\pi$ should be the Artin ${\rm L}$-function of the corresponding $\rho$, with matching conductors, the weights of $\pi$ corresponding to the Hodge-Tate weights $\rho$, and with ${\rm w}(\pi)$ equal to the Deligne weight of $\rho$.  For example, $\rho$  is of finite image (Artin type), if and only if the weights of $\pi$ are all $0$. \ps

A general finiteness result of Harish-Chandra \cite{harishchandra} ensures that for any $n\geq 1$, there are only finitely many cuspidal algebraic $\pi$ of ${\rm GL}_n$ over $E$, with given weights and conductor. Nevertheless, beyond a few situations in which discrete series representations of are involved, it seems hard to say much about the number of those representations, even in simple specific cases such as dimension $n=3$, $E=\Q$ and conductor $1$. Our main result in this paper is a finiteness theorem which is {\it uniform in  the dimension $n$}. It is both easier to state, and stronger, in the case $E$ is the field of rational numbers. \ps


\begin{thmintro}\label{thmch} Let $N\geq 1$ be an integer. There are only finitely many cuspidal algebraic automorphic representations $\pi$ of $\GL_n$ over $\Q$, with $n$ varying, whose conductor is $N$, such that the weights of $\pi$ are in the interval $\{0,\dots,23\}$. \end{thmintro}
\noindent 
\ps

This finiteness statement is surprising to us and we are not aware of any {\it a priori} philosophical reason, either automorphic, nor Galois theoretic or motivic, why it should hold. It asserts in particular that $N$ being given, then for $n$ big enough there is no cuspidal automorphic representation as in the statement. This consequence is even equivalent to the full statement, by the aforementioned result of Harish-Chandra. Assuming the yoga of Fontaine-Mazur and Langlands, the finiteness results discussed in the beginning of the introduction  
correspond to the special cases of Theorem \ref{thmch} with $\{0,\dots,23\}$ replaced with $\{0,1\}$. \ps

Before discussing the proof of Theorem \ref{thmch}, we state a generalization which holds for number fields $E$ whose root-discriminant ${\rm r}_E = |{\rm disc}\, E|^\frac{1}{[E:\Q]}$ is sufficiently small. For any integer $w\geq 0$, consider the symmetric matrix ${\rm M}(w)= (\,{\rm log}\, \pi\, -\, \psi(\,\frac{1+|i-j|}{2}\,)\,)_{0 \leq i,j \leq w}$, where $\psi(s)=\frac{\Gamma'(s)}{\Gamma(s)}$ denotes the classical digamma function. We also denote by ${\rm c}(M)$ denotes the sum of all the coefficients of the matrix $M$, and we define real numbers ${\rm t}(w)$ and ${\rm r}(w)$ by the formulas
$${\rm t}(w) = 1/{\rm c}({\rm M}(w)^{-1}) \,\,\,\,\,\,{\rm and} \,\,\,\,\,\,{\rm r}(w) = {\rm exp}({\rm t}(w))$$
(${\rm t}(w)$ is actually well-defined). We also set\footnote{Recall the simple formula $\psi(1+w)=-\gamma+\sum_{1 \leq k\leq w} \frac{1}{k}$ for $w\geq 0$ an integer, where $\gamma$ is the {\it Euler-Mascheroni constant}.} 
${\rm r}^\ast(w)=8\pi e^{-\psi(1+w)}$. We will show that ${\rm r}(w)$ and ${\rm r}^\ast(w)$ are nonincreasing functions of $w$, with ${\rm r}^\ast(w) > {\rm r}(w)$ (see \S\ref{pfmainthm}). Our main result is the following:
\ps
\begin{thmintro}\label{mainthm} Let $w\geq 0$ be an integer, $E$ a number field whose root-discriminant ${\rm r}_E$ satisfies ${\rm r}_E < {\rm r}(w)$, and $\mathcal{N} \subset \mathcal{O}_E$ an ideal.  There are only finitely many cuspidal algebraic automorphic representations $\pi$ of $\GL_n$ over $E$, with $n$ varying, whose conductor is $\mathcal{N}$ and whose weights are in the intervall $\{0,\dots,w\}$.\ps
Moreover, under a certain version ${\rm (GRH)}$ of the Generalized Riemann Hypothesis, the same result holds if we replace ${\rm r}_E < {\rm r}(w)$ by the weaker condition ${\rm r}_E < {\rm r}^\ast(w)$. 
\end{thmintro}
\ps\ps

\noindent  The table below gives the relevant\footnote{Classical estimates of Minkowski show ${\rm r}_E > 2$ for $[E:\Q]>2$; this explains why the useless informations concerning $15 \leq w \leq 22$ are omitted in Table \ref{tablerw}.}  values of ${\rm r}(w)$ and ${\rm r}^\ast(w)$:
\ps
\begin{table}[!h]
\renewcommand{\arraystretch}{1.5}
{\scriptsize
{
\begin{tabular}{c||c|c|c|c|c|c|c|c|c}
 $w$  & $0$ & $1$ & $2$ & $3$ & $4$ & $5$ & $6$ & $7$ & $8$ \\ 
\hline ${\rm r}(w)$ & $22.3816$ & $11.1908$ & $7.5690$ & $5.7456$ & $4.6401$ & $3.8959$ & $3.3597$ & $2.9546$ & $2.6375$ \\
\hline ${\rm r}^\ast(w)$ & $44.7632$ & $16.4675$ & $9.9880$ & $7.1567$ & $5.5737$ & $4.5633$ & $3.8628$ & $3.3486$ & $2.9551$ \\
\hline $w$ & $9$ & $10$ & $11$ & $12$ & $13$ & $14$ & $\cdots$ & $23$ & $24$  \\
\hline ${\rm r}(w)$ & $2.3824$ & $2.1726$ & $1.9971$ & $1.8480$ & $1.7197$ & $1.6082$ & $\cdots$ & $1.0167$ & $0.9768$ \\
\hline ${\rm r}^\ast(w)$ & $2.6443$ & $2.3927$ & $2.1848$ & $2.0101$ & $1.8613$ & $1.7329$ & $\cdots$ & $1.0694$ & $1.0258$ \\
\end{tabular}\ps\ps
\caption{{\small Values of ${\rm r}(w)$ and ${\rm r}^\ast(w)$ up to $10^{-4}$.}}
\label{tablerw}
}
}\end{table}

\noindent For instance, we have ${\rm r}(23) > {\rm r}_\Q=1$ and Theorem \ref{mainthm} does imply Theorem \ref{thmch}. We also have ${\rm r}^\ast(24)>1$ so Theorem \ref{mainthm} asserts that we may replace $23$ by $24$ in the statement of Theorem \ref{mainthm} if we assume ${\rm (GRH)}$; note however ${\rm r}^\ast(25)<1$. For $E=\Q(\sqrt{-3})$, we have ${\rm r}_{E} = \sqrt{3} \approx 1.7321$ up to $10^{-4}$, so the table shows that the finiteness statement holds with $w=12$ (resp. $w=14$ under ${\rm (GRH)}$).We also mention ${\rm r}(1)=2\pi e^{\gamma}$.

\subsection{Ideas of the proofs and organization of the paper} The proofs of Theorems \ref{thmch} and \ref{mainthm} are of analytic flavor, and very much inspired by the seminal works of Stark, Odlyzko and Serre on discriminant bounds \cite{poitou,odlyzko}, continued by Mestre in \cite{mestre}. \ps
In order to prove, say, Theorem \ref{thmch}, a natural idea is to study the {\it explicit formula} ``\`a la Weil''  for the Godement-Jacquet ${\rm L}$-function of a cuspidal algebraic $\pi$ of ${\rm GL}_n$ over $\Q$ of conductor $N$. Using Odlyzko's {\it test functions}, this gives a non trivial lower bound of the form ${\rm log}\, N \,\geq \alpha \, n$, for some explicit real number $\alpha>0$, provided we have ${\rm w}(\pi) \leq 10$. This does proves Theorem \ref{thmch} (and much more) in this range, but we seem to be stuck for higher ${\rm w}(\pi)$ and still quite far from the high bound $23$ announced in Theorem \ref{thmch}. The situation is even far worse for Theorem \ref{mainthm} when the number field $E$ is different from $\Q$: this method gives nothing for ${\rm r}_E> 2.67$ (see \S\ref{retoursurgj}). \ps
	
	In this work we examine the explicit formula associated to the Rankin-Selberg ${\rm L}$-function of $\pi \times \pi^\vee$, an idea which goes back at least to Serre in the Artin case \cite[\S 8]{poitou} and which has also been used previously by others in the past in related contexts (e.g. by Miller in \cite{miller}). This formula depends on the choice of a {\it test function} $F : \R \rightarrow \R$ (\S \ref{testfunctions}), and it is well-known that under a suitable positivity assumption on both $F$ and its Fourier transform, but without any assumption on the support of $F$, it leads to a simple yet nontrivial inequality depending on $F$. Combined with the general bounds for conductors of pairs due to Bushnell and Henniart \cite{bh}, this inequality may be written  (see \S \ref{sectrs}) \begin{equation} \label{basinegintro} {\rm q}_F^\R (V) \leq \widehat{F}(\frac{i}{4\pi}) \,+\, F(0) \, \dim V \, {\rm log} \, N.\end{equation}
Here, ${\rm q}_F^\R$ is a certain real-valued quadratic form depending only on $F$, defined on the Grothendieck ring ${\rm K}_\R$ of the maximal compact quotient of the Weil group of $\R$ (an extension of $\Z/2\Z$ by the circle), and $V$ is the class in ${\rm K}_\R$ of the Langlands parameter of {\small $\pi_\infty \otimes |\det|^{-\frac{{\rm w}(\pi)}{2}}$}. In particular, we have $\dim V=n$. \ps
 
 The basic idea of our proof is as follows. If the representation $\pi$ has nonnegative weights, and motivic weight $w$, then $V$ lies in a specific finite rank sublattice ${\rm K}_\R^{\leq w}$ in ${\rm K}_\R$. Assume we may find an allowed test function $F$ such that ${\rm q}_F^\R$ is positive definite on ${\rm K}_\R^{\leq w}$. As the set of points of an Euclidean lattice whose squared norm is less than a given affine function is finite, the inequality \eqref{basinegintro} implies that $V$ lies in a finite subset of ${\rm K}_\R^{\leq w}$, which in turn implies that $\dim V$ is bounded and the finiteness of the possible representations $\pi$ by Harish-Chandra. \ps
 
 	This explains why an important part of this work is devoted to studying the family of quadratic forms ${\rm q}^\R_F$, as well as a related family of quadratic forms ${\rm q}_F^\C$ on the Grothendieck ring ${\rm K}_\C$ of the unit circle. Note that when $\pi$ is a cuspidal algebraic automorphic representation of ${\rm GL}_n$ over a general number field $E$, a variant of \eqref{basinegintro} holds, but with two main differences. There is first an additional term on the right of the form $\frac{F(0)}{2}\, n^2\,|{\rm disc}\, E|$, explained by Formula \eqref{formcondabs}. Moreover, the left-hand side becomes a sum of terms, indexed by the Archimedean places $v$ of $E$, of the form ${\rm q}_F^{E_v}(V_v)$, where $V_v$ is the class in ${\rm K}_{E_v}$ of {\small $\pi_v \otimes |\det|_v^{-\frac{{\rm w}(\pi)}{2}}$}. The same strategy as above can be applied, provided we can study the signature of the quadratic form ${\rm q}_F^L\, -\, \frac{F(0)}{2} \,[L:\R]\,{\rm log}\, {\rm r}_E \,\dim^2$ on the lattices ${\rm K}_L^{\leq w}$ for $L=\R$ or $\C$.  Actually, as explained in Proposition \ref{wrwc}, the study of the real case may be entirely deduced from the complex one, when $F$ is nonnegative; as the ring ${\rm K}_\C$ is a little simpler than ${\rm K}_\R$ we do focus on the complex case. \ps
	
		The whole of section \S \ref{sect1} is devoted to studying, for any $t \in \R$, the quadratic forms ${\rm q}_F^\C - t \dim^2$ on ${\rm K}_\C$ and its natural filtration by the ${\rm K}_\C^{\leq w}$ for $w \geq 0$, for a natural class of functions $F$ containing all test functions. We work in a context independent of the discussion above and use mostly here elementary harmonic analysis on the circle. Denote\footnote{We actually have ${\rm K}_L^{\leq w} \subset {\rm K}_L^{\leq w+2}$ for any $w\geq 0$, and ${\rm K}_L^{+}$ (resp. ${\rm K}_L^{-}$) is the union of the ${\rm K}_L^{\leq w}$ with $w$ even (resp. odd).} by ${\rm K}_L^{\pm} \subset {\rm K}_L$ the $\pm$-eigenspace of the central element $-1$ in ${\rm W}_L$. A first striking result is that for any nonzero nonnegative $F$, and for $L=\R$ and $\C$, the quadratic form ${\rm q}_F^L$ is positive definite on the hyperplane of ${\rm K}_L^{\pm}$ defined by $F(0) \dim = 0$  (Proposition \ref{posdefTH}); it has signature $(\infty,1)$ on the whole of ${\rm K}_L^{\pm}$ if we have $F(0) \neq 0$. This statement contains a collection of elementary yet nontrivial inequalities, such as the following one (for a particular $F$): for any real numbers $x_0,\dots,x_n$ with $\sum_{i=0}^n x_i=0$ we have $\sum_{0 \leq i,j \leq n} \, x_i x_j \,{\rm log}\, (1+|i-j|) \leq 0$. \ps
		
		Two specific functions $F$ play an important role, namely $F=1$ and $F(t)={\rm cosh}(\frac{t}{2})^{-1}$. They are not test functions, but the respective limits of the two families of Odlyzko's test functions ${\rm G}_\lambda$ and ${\rm F}_\lambda$, with $\lambda \rightarrow \infty$ (see \S \ref{basicinequality}, the function ${\rm G}_\lambda$ is allowed in \eqref{basinegintro} only under {\rm (GRH)}). The story of the form ${\rm q}_1^\C$ is especially beautiful. It is luckily related to the family of Legendre orthogonal polynomials, which allows us to show that ${\rm q}_1^\C- t \dim^2$ is positive definite on ${\rm K}_\C^{\leq w}$ if, and only if, we have $\psi(1+w)<{\rm log}\, 8\pi -t$ (see \S \ref{casHgrh}). As we have the equivalence
$$1+\frac{1}{2}+\frac{1}{3}+\dots+\frac{1}{w} < {\rm log}\, 8\pi +\gamma \,\,\, \Leftrightarrow \,\,\, w \leq 24,$$ this shows that ${\rm q}_1^L$ is positive definite on ${\rm K}_L^{\leq w}$ if and only if we have $w \leq 24$ ($L=\R$ or $\C$). Under {\rm (GRH)}, these facts imply Theorem \ref{thmch} with $24$ replacing $23$, and eventually Theorem \ref{mainthm}, and gives as well a conceptual explanation of the absolute bound $w\leq 24$ in all of our results, without relying on any sophisticated numerical computation. The study of ${\rm q}_F$ with $F(t)={\rm cosh}(\frac{t}{2})^{-1}$ is more delicate (see \S \ref{casHsansgrh}), and does rely on a few numerical computations to obtain Table \ref{tablerw}, but they are very easy to justify. This eventually leads to the proof of both Theorems. Let us mention that ${\rm M}(w)$ is a Gram matrix of ${\rm q}_F$ on ${\rm K}_\C^{\leq w}$. \ps

\subsection{Other results, perspectives and applications} It is a very interesting open question to us whether Theorem \ref{thmch} should hold in motivic weight $25$ or higher. The reader may object that in our proof we used very specific test functions, but perhaps other choices would lead to stronger results. We show in \S \ref{chaptopt} that this is not the case, so that our results are optimal in this sense (for any $E$). For instance, we prove that for any test function $F$ satisfying $F \geq 0$ and $\widehat{F} \geq 0$, then ${\rm q}_F$ is not positive definite on ${\rm K}_\R^{\leq 25}$. Worse, there are concrete effective elements in ${\rm K}_\R^{\leq 25}$ on which ${\rm q}_F$ has a negative value for all nonzero such $F$. For example, we cannot rule out the possibility that for some integer $N\geq 1$, there exist infinitely many integers $m\geq 1$ with a cuspidal algebraic automorphic representation of $\GL_{32 m}$ over $\Q$, of conductor $N$, whose weights are $1$, $2$, $\dots$, $24$ with multiplicity $m$, and both $0$ and $25$ with multiplicity $4m$. As a complement, we also show in \S \ref{chaptopt} that without assuming ${\rm (GRH)}$, our method cannot prove the statements that we do prove assuming ${\rm (GRH)}$ !\ps
   
Theorems~\ref{thmch} and \ref{mainthm} naturally raise the question of determing the set ${\rm S}(E,\mathcal{N},w)$ of cuspidal algebraic representations $\pi$ of ${\rm GL}_n$ over $E$ of given conductor $\mathcal{N}$ and weights in $\{0,\dots,w\}$, when we have ${\rm r}_E < {\rm r}(w)$. As our proof is effective in several respects (see \S \ref{chaptopt}), especially for a trivial $\mathcal{N}$, it provides a starting point for this question. This is the first step in the determination of ${\rm S}(\Q,(1),22)$ by J. Lannes and the author in \cite[Chap. IX]{CL}. This set has exactly $11$ elements up to twist, a result which has interesting applications to automorphic forms on classical groups over $\Z$ as shown in \cite{CL} (such as, a new proof that there are exactly $24$ isometry classes of even unimodular lattices of rank $24$). This is a very specific case and there is still much to explore in general. In this direction, we mention some forthcoming work by O. Ta\"ibi and the author about ${\rm S}(\Q,(1),w)$ for $w=23$ and $24$, as well as the Ph. D. thesis of G. Lachauss\'ee for some results about ${\rm S}(\Q,(N),w)$ with $N>1$. In another direction, it should be possible to prove variants of our main results for non necessarily algebraic automorphic representations. \ps\ps

{\small \noindent {\sc Acknowledgement.} This work was started during the joint work \cite{CL} with Jean Lannes, and it is a pleasure to thank him here, as well as Guillaume Lachauss\'ee, Patrick G\'erard and Olivier Ta\"ibi, for several useful discussions.\ps\ps}

{\small
\tableofcontents
}

\section{Some quadratic forms}\label{sect1}

\subsection{} \label{somequad} We fix a Lebesgue integrable function $H : [0,+\infty[ \rightarrow \R$ such that the map $t \mapsto \frac{H(t)-H(0)}{t}$ is bounded on $]0,\epsilon[$ for some $\epsilon>0$. The sum
\begin{equation}\label{defpsih} \psi_H(s) \,=\, \int_0^\infty \left( H(0) \frac{e^{-t}}{t} - H(t) \frac{e^{-st}}{1-e^{-t}} \right)\,{\rm d}t, \end{equation}
is then absolutely convergent for all $s \in \C$ with ${\rm Re} \,s \geq 0$, and defines a holomorphic function on ${\rm Re}\,s\,>0$ that is continuous on ${\rm Re}\, s \geq 0$. For example, for $a>0$ and $H(t)=e^{-at}$ we have $\psi_H(s)=\psi(a+s)$ where $\psi(s)=\frac{\Gamma'(s)}{\Gamma(s)}$ is the classical {\it digamma} function, and for $H(t) = \frac{1-e^{-t}}{t}e^{-at}$ we have the well-known identity $\psi_H(s)= {\rm log}\, (a+s)$. 
For $n\geq 0$ an integer, we consider the quadratic form ${\rm q}_n^H : \R^{n+1} \rightarrow \R$ defined by 
\begin{equation} \label{defqnh} {\rm q}_n^H(x_0,\dots,x_n)=\sum_{0 \leq i,j \leq n }\, x_i\, x_j \,\psi_H(\,|i-j|\,). \end{equation} 
\begin{prop}\label{posdefTH} Assume $H \geq 0$ and that $H$ is nonzero in ${\rm L}^1(]0,+\infty[)$. Then for all integers $n\geq 0$, the quadratic form ${\rm q}_n^H$ is negative definite on the subspace $H(0) \, \sum_i x_i=0$ of $\R^{n+1}$. \end{prop}
\ps\ps
In other words, ${\rm q}_n^H$ is negative definite on the hyperplane $\sum_i x_i=0$, and even on the whole of $\R^{n+1}$ if we have $H(0) = 0$. In order to prove this proposition we introduce the distribution on $\R/\Z$ defined by
\begin{equation} \label{defTH} {\rm T}_H  \,=\, \sum_{n \in \Z} \, \, \, \psi_H(\,|n|\,) \,\,e^{2i\pi nx}.\end{equation}
This is an even distribution, which is of order $\leq 2$ as we have 
\begin{equation} \label{asym} \psi_H(1+s) - \,H(0) \,{\rm log} (1+s) \rightarrow 0,\, \, \, \, \, \, \, \, \, \,  {\rm Re}\, s \rightarrow +\infty \end{equation}
(reduce to the case $H(0)=0$ by replacing $H(t)$ with the function $H(t)-H(0) \frac{1-e^{-t}}{t}e^{-t}$). It is of interest because of the obvious identity 
\begin{equation} \label{relbase} {\rm q}_n^H(x_0,\dots,x_n)\,=\, \langle {\rm T}_H , | \sum_{k=0}^n \,x_k \,e^{2i\pi kx}\,|^2\,\rangle, \, \, {\rm with}\, \,  x_0,\dots,x_n \in \R.\end{equation}

\begin{lemme} \label{identTH} For any smooth $1$-periodic function $\varphi$ on $\R$ we have 
{\scriptsize 
$$\langle {\rm T}_H,\varphi\rangle \,= \, \, \int_{]0,\infty[\times ]0,1[} \frac{ 1-e^{-2 t}}{|1-e^{-t+2i\pi x}|^2}\,\left(\varphi(0) H(0) \frac{e^{-t}}{t} - \varphi(x)H(t)\frac{1}{1-e^{-t}}\right) \, {\rm d}t\,{\rm d}x.$$
}
\end{lemme} 

 Let us show first that Lemma~\ref{identTH} implies Proposition~\ref{posdefTH}. Fix some elements $x_0,\dots,x_n \in \R$ with $H(0)(\sum_k x_k)=0$ and consider $\varphi(x) = |\sum_{k=0}^n \,x_k\, e^{2i k\pi x}|^2$. We have $\varphi \geq 0$ and $H(0) \varphi(0)=0$. Lemma~\ref{identTH} thus shows $\langle {\rm T}_H , \varphi \rangle  \leq 0$ for $H \geq 0$. This inequality is strict if furthermore $\varphi \neq 0$ and $H$ is nonzero in ${\rm L}^1(]0,+\infty[)$, as $\varphi$ has only finitely many zeros in $\R/\Z$. This proves Proposition \ref{posdefTH} by Formula \eqref{relbase}. $\square$

\begin{pf} (of Lemma \ref{identTH}) Let $\psi$ be a Lebesgue integrable function on $[0,+\infty[ \times [-\frac{1}{2},\frac{1}{2}]$. Assume we have $|\psi(t,x)|\leq a |t| + b |x|$ for all $(t,x)$ in some neighborhood of $(0,0)$, with fixed $a,b \in \R$. Then we have 
\begin{equation}\label{converg} \int_{[0,+\infty[ \times [-\frac{1}{2},\frac{1}{2}]} \frac{|\psi(t,x)|}{|1-e^{-t+2i\pi x}|^2} \, {\rm d}t{\rm d}x < \infty.\end{equation}
Indeed, this integrability is obvious outside any neighborhood of $(0,0)$ in $[0,+\infty[ \times [-\frac{1}{2},\frac{1}{2}]$. But we have the equivalence $|1-e^{-t+2i\pi x}|^2 \sim t^2+ (2\pi x)^2$ for $(t,x) \rightarrow (0,0)$, so \eqref{converg} follows from the assumption on $\psi$ and the integrability of $(u,v) \mapsto \frac{u}{u^2+v^2}$ on the disc $u^2+v^2<1$. The bound \eqref{converg} applies in particular to $\psi(t,x)=\varphi(0) H(0) \frac{e^{-t}(1-e^{-t})}{t} - \varphi(x)H(t)$, as $t \mapsto \frac{H(t)-H(0)}{t}$ bounded in a neighborhood of $0$ by assumption: the integral of the statement is absolutely convergent.\ps
 
Let us now view ${\rm T}_H$ as a $2$-periodic distribution. Multiplying ${\rm T}_H$ by $e^{-i\pi x}-e^{i\pi x}$ we have the identity
\begin{equation} \label{eqan1} - 2 i\, {\rm sin} \pi x\, {\rm T}_H\, =\, \sum_{n \in \Z} \, a_n \, e^{i \pi n x}\end{equation}
with $a_{n}=- a_{-n}$, $a_n=0$ for $n$ even, and for $n\geq 1$ odd 
$$a_n = \psi_H(\frac{n+1}{2}) - \psi_H(\frac{n-1}{2}) =  \int_0^\infty H(t) e^{- \frac{n-1}{2} t} {\rm d}t.$$
We then proceed to a geometric summation. For $t>0$ we have 
$$\Im\, \sum_{k=0}^\infty e^{-kt}e^{i(2k+1)\pi x} =  \Im\, \frac{e^{i\pi x}}{1- e^{-t + 2 i \pi x}} = {\rm sin}\pi x  \,\frac{ \, 1+e^{-t} }{|1-e^{-t+2i\pi x}|^2},$$
as well as $|\Im\, e^{i\pi x}\,\sum_{k=0}^{n} e^{-kt} e^{2i \pi kx}| \leq 4 |1-e^{-t+2i\pi x}|^{-2}$, for all $x \in \R$ and $n\geq 0$. The bound \eqref{converg} and the dominated convergence theorem imply
$$ \langle  \,{\rm sin}\, \pi x\, {\rm T}_H, \varphi \rangle \,=\, -\frac{1}{2}\, \int_{]0,\infty[\times ]0,2[} H(t)\,\varphi(x)\, {\rm sin}\, \pi x\,  \frac{1+e^{-t}}{|1-e^{-t+2i\pi x}|^2}\,{\rm d}t{\rm d}x$$
for any smooth function $\varphi$ on $\R/2\Z$ such that $\varphi(0)=\varphi(1)=0$, {\it i.e.} which is divisible by ${\rm sin} \pi x$. This proves the equality of the lemma in the case $\varphi$ is divisible by $({\rm sin}\, \pi x)^2$, {\it i.e.} with $\varphi(0)=\varphi'(0)=0$. As both sides of this equality are linear in $\varphi$, it only remains to prove it in the cases $\varphi(x)=1$ and $\varphi(x)={\rm sin}\, 2 \pi x$. In the latter case, both sides are $0$ for a parity reason. In the case $\varphi=1$, we have 
$\langle {\rm T}_H, 1\rangle = \psi_H(0) = \int_0^\infty \left(H(0) \frac{e^t}{t} - H(t)\frac{1}{1-e^{-t}}\right){\rm d}t$,
by the definitions of ${\rm T}_H$ and of $\psi_H$. We conclude by the simple identity $\int_0^1 \frac{1-e^{-2t}}{|1-e^{-t+2i\pi x}|^2} {\rm d}x = 1$, which holds for all $t>0$. \end{pf}


\begin{example} \label{examplelog}{\rm Consider $H(t)=e^{-t}(1-e^{-t})/t$, in which case we have $\psi_H(s)={\rm log}(1+s)$. Proposition \ref{posdefTH} applies and shows that the signature of $(\,{\rm log} \,(1+|i-j|)\,)_{0\leq i,j\leq n}$ is $(1,n)$ for all $n>0$ (surprisingly!). We also have  $\psi_{tH}(s)=- \frac{\partial}{\partial s}\psi_H(s) = -\frac{1}{1+s}$. By Proposition \ref{posdefTH} again, the matrix $(\frac{1}{1+|i-j|})_{0\leq i,j\leq n}$ is positive definite for all $n\geq 0$.} \end{example}

\subsection{The full signature} Let $H$ be as in \S \ref{somequad}; we now assume $H\geq 0$ and that $H$ is nonzero in ${\rm L}^1([0,+\infty[)$. Fix an integer $n\geq 0$ and denote by $\phi_n$ the linear form on $\R^{n+1}$ defined by 
\begin{equation}\label{eefphin} \phi_n(x_0,\dots,x_n)=\sum_{0 \leq k \leq n} x_k. \end{equation} We are interested in the signature of the quadratic form\footnote{Note that this form corresponds to the $1$-periodic distribution ${\rm T}_H - t \delta_0$ via an identity similar to \eqref{relbase}.} $t \phi_n^2 - {\rm q}_n^H $ on $\R^{n+1}$, for $t \in \R$. On the hyperplane ${\rm ker}\, \phi_n$, this quadratic form coincides with $-{\rm q}_n^H$, hence is positive definite by Proposition \ref{posdefTH}. In particular, there is a unique element ${\rm v}_n^H$ in $\R^{n+1}$ which is orthogonal to ${\rm ker}\, \phi_n$ with respect to ${\rm q}_n^H$, and with $\phi_n({\rm v}_n^H)=1$. Set
\begin{equation} \label{defth} {\rm t}_n^H = {\rm q}_n^H({\rm v}_n^H).\end{equation}
\begin{prop}\label{crittnh} Assume $H\geq 0$ and that $H$ is nonzero in ${\rm L}^1(]0,+\infty[)$. Let $n\geq 0$ be an integer and $t \in \R$. \begin{itemize} \ps
\item[(i)] The quadratic form $t \phi_n^2- {\rm q}_{n}^H$ is positive definite for $t>{\rm t}_n^H$, degenerated for $t={\rm t}_n^H$, and of signature $(n,1)$ for $t<{\rm t}_n^H$.\ps\ps
\item[(ii)] We have the inequalities $\frac{1}{2}(\psi_H(0)+\psi_H(n)) \,\leq \, {\rm t}_n^H \leq {\rm t}_{n+1}^H$. \ps\ps
\item[(iii)] If ${\rm t}_n^H \neq t$, then $1/(t-{\rm t}_n^H)$ is the sum of the coefficients of the inverse of the matrix $(t-\psi_H(\,|i-j|\,))_{0 \leq i,j \leq n}.$\end{itemize}
\end{prop}

\begin{pf} As the decomposition $\R^{n+1} \,=\, {\rm ker}\, \phi_n \,\oplus \,\R \,{\rm v}_n^H$ is orthogonal for both ${\rm q}_n^H$ and $\phi_n^2$, assertion (i) is immediate. We now prove (ii). Fix $u \in \R$ with $u>{\rm t}_{n+1}^H$; then $u \phi_{n+1}^2- {\rm q}_{n+1}^H$ is positive definite on $\R^{n+2}$ by assertion (i). The restriction of this quadratic form to $\R^{n+1}$, viewed as the subspace of $\R^{n+2}$ with last coordinate $0$, is $u\phi_n^2-{\rm q}_{n}^H$, and negative definite, so we have $u> {\rm t}_n^H$ by (i) again. This shows ${\rm t}_n^H \leq {\rm t}_{n+1}^H$. For the other inequality, use ${\rm q}_n^H(x) \,\leq \,{\rm t}_n^H \,\phi_n(x)^2$ for $x=(1,0,\dots,0,1)$.
\par
 We now prove (iii). The Gram matrix ${\rm J}_n$ of $(x,y) \mapsto \phi_n(x)\phi_n(y)$ in the canonical basis of $\R^{n+1}$ has all its coefficients equal to $1$. For $t \neq {\rm t}_n^H$, the form $t \phi_n^2-{\rm q}_n^H$ is nondegenerate by (i), so its Gram matrix $\gamma_n$ in the canonical basis of $\R^{n+1}$, which is the matrix in assertion (iii), is invertible. By (i) and the analysis before the proposition, the determinant of $\gamma_n-u {\rm J}_n$ vanishes exactly for $u=t-{\rm t}_n^H$, and is linear in $u$, hence equal to $(\det \gamma_n)\, (1-u \,{\rm Tr}(\gamma_n^{-1} {\rm J}_n)\,)$.
\end{pf}


\begin{prop} \label{inegqh1h2}Let $H_i$ be as in \S \ref{somequad} for $i=1,2$, and assume $H_1 \leq H_2$ and $H_1(0)=H_2(0)$. Then we have ${\rm q}_n^{H_2} \leq {\rm q}_n^{H_1}$ for all $n\geq 0$. If furthermore we have $H_i \geq 0$ and $H_i$ nonzero in ${\rm L}^1(]0,+\infty[)$ for $i=1,2$, then we have ${\rm t}_n^{H_2} \leq {\rm t}_n^{H_1}$ for all $n\geq 0$.
\end{prop}
\ps
\begin{pf} Fix an integer $n\geq 0$. We have ${\rm q}_n^{H_2}= {\rm q}_n^{H_1} + {\rm q}_n^{H_3}$ with $H_3=H_2-H_1$. The assumptions $H_3(0)=0$ and $H_3\geq 0$, as well as Proposition \ref{posdefTH}, imply ${\rm q}_n^{H_3} \leq 0$ (of course, we have ${\rm q}_n^{H_3}=0$ if $H_3$ is almost everywhere $0$). This shows the first assertion. This assertion implies  ${\rm t}_n^{H_2} \phi_n^2- {\rm q}_n^{H_1} \leq {\rm t}_n^{H_2} \phi_n^2- {\rm q}_n^{H_2}$, hence the second assertion by Proposition \ref{crittnh} (i).\end{pf}

\subsection{} Our final aim in this section is to study the element ${\rm v}_n^H$ and the number ${\rm t}_n^H$ in the two specific cases that we shall need later, namely with $H(t)=e^{-\frac{t}{2}}$ or $H(t)\,=\,e^{-\frac{t}{2}}\,{\rm cosh}(\frac{t}{2})^{-1}$. We first continue the analysis made in \S \ref{somequad}, assuming $H\geq 0$, and $H \neq 0$ in ${\rm L}^1(]0,+\infty[)$. \ps

Let $\mathcal{P}$ be the $\R$-algebra of $2$-periodic trigonometric polynomials with real coefficients: we have $\mathcal{P}=\R[e^{i\pi x},e^{-i\pi x}]$ and $\mathcal{P} = \cup_{n\geq 0} \mathcal{P}_n$, where $\mathcal{P}_n$ is the subspace of homogeneous polynomials of degree $n$ in $e^{i\pi x}$ and $e^{-i\pi x}$. We have $\dim \mathcal{P}_n=n+1$, $\mathcal{P}_n \subset \mathcal{P}_{n+2}$, and $\mathcal{P}_n\mathcal{P}_m \subset \mathcal{P}_{n+m}$. The involution of $\mathcal{P}$ sending $P(x)$ to its complex conjugate $\overline{P}(x)=P(-x)$ preserves each $\mathcal{P}_n$. We denote by $\mathcal{P}^+ \subset \mathcal{P}$ the subalgebra of real valued elements; the theory of Tchebychev polynomials shows $\mathcal{P}^+=\R[{\rm cos}\,\pi x]$.\ps

We may and do view ${\rm T}_H$ as a $2$-periodic distribution, hence as a linear form on $\mathcal{P}$. As ${\rm T}_H$ is real, the formula $(P, Q) \mapsto \langle {\rm T}_H, P \overline{Q} \rangle$ defines a symmetric bilinear form on $\mathcal{P}$. Under the isomorphism $(x_0,\dots,x_n) \mapsto e^{-i\pi n x} \sum_{k=0}^n x_k e^{2i\pi k x}$ between $\R^{n+1}$ and $\mathcal{P}_n$, the quadratic form ${\rm q}_n^H$ and the linear form $\phi_n$ correspond respectively to $P \mapsto \langle {\rm T}_H, P \overline{P} \rangle$ and $P \mapsto P(0)$. We denote by ${\rm P}_n^H$ the element of $\mathcal{P}_n$ corresponding to ${\rm v}_n^H \in \R^{n+1}$ under the isomorphism above. We have ${\rm P}_n^H(0)=1$ by definition. \ps\ps

 
\begin{prop}\label{orthpol} Assume $H\geq 0$ and that $H$ is nonzero in ${\rm L}^1(]0,+\infty[)$. For any $n\geq 0$, the trigonometric polynomials ${\rm P}_0^H, {\rm P}_1^H, \dots,{\rm P}_n^H$ are real and pairwise orthogonal with respect to the inner product on $\mathcal{P}^+$ 
\begin{equation} \label{innerprodp} (P,Q) \mapsto -\langle ({\rm sin} \,\pi x)^2 {\rm T}_H, PQ \rangle.\end{equation} Moreover, ${\rm P}_n^H$ has ``degree'' $n$ and we have ${\rm t}_n^H = \langle {\rm T}_H\, ,\,e^{i\pi n x} {\rm P}_n^H \rangle$.
\end{prop}

\begin{pf} Note that Formula \eqref{innerprodp} defines an inner product on $\mathcal{P}^+$ by Lemma \ref{identTH} and our assumption on $H$. By definition, ${\rm P}_n^H$ is the unique element of $\mathcal{P}_n$ with ${\rm P}_n^H(0)=1$ and $ \langle{\rm T}_H, {\rm P}_n^H \overline{Q} \rangle=0$ for all $Q$ in $\mathcal{P}_n$ with $Q(0)=0$. This uniqueness property implies ${\rm P}_n^H = \overline{{\rm P}^H_n}$. The element $f=e^{i\pi x}-e^{-i\pi x}\,=\,2i \,{\rm sin} \pi x$ is in $\mathcal{P}_1$, so we have $f^2 \mathcal{P}_m \subset \mathcal{P}_{m+2}$, and $f$ vanishes at $0$. This shows the orthogonality \begin{equation} \label{orthopnh}  - \langle\, ({\rm sin} \pi x)^2\,{\rm T}_H, {\rm P}_n^H {\rm P}_m^H\,\rangle\,=\,\frac{1}{4}\, \langle\, {\rm T}_H, {\rm P}_n^H   \,f^2\,{\rm P}_m^H\,\rangle\,=\,0 \end{equation} for $0 \leq m<n$ and $m \equiv n \bmod 2$. On the other hand, we have $\langle {\rm T}_H, P \rangle = 0$ for every $P$ in $\mathcal{P}_n$ with $n$ odd, by the $1$-periodicity of ${\rm T}_H$ and the identity $P(x+1)=-P(x)$. This shows \eqref{orthopnh} for all $m \not \equiv n \bmod 2$, and concludes the proof of the first assertion. The fact that ${\rm P}_n^H$ has degree exactly $n$, {\it i.e.} does not belong to $\mathcal{P}_m$ for $m<n$, follows by induction from this assertion. Last but not least, we have ${\rm t}_n^H = \langle {\rm T}_H, {\rm P}_n^H \overline{{\rm P}_n^H}\rangle$ by definition. As ${\rm P}_n^H-e^{-i \pi n x}$ is an element of $\mathcal{P}_n$ vanishing at $0$, we have $\langle {\rm T}_H, {\rm P}_n^H \overline{{\rm P}_n^H - e^{-i \pi n x}} \rangle =0$, and we are done. 
\end{pf}

\subsection{The case $H(t)=e^{-\frac{t}{2}}$}\label{casHgrh} In this case, we have $\psi_H(s)=\psi(\frac{1}{2}+s)$. 

\begin{lemme}\label{sinTHgrh} For $H(t)=e^{-\frac{t}{2}}$, we have $-({\rm sin}\, \pi x)^2 \,{\rm T}_H \,=\, \frac{\pi}{2}|{\rm sin} \,\pi x|$. 
\end{lemme}

\begin{pf} We apply Formula \eqref{eqan1}. In our specific case, we have $a_n=\psi(1+\frac{n}{2})-\psi(\frac{n}{2})=\frac{2}{n}$ for $n\geq 1$ odd, by the identity $\psi(1+x)=\frac{1}{x}+\psi(x)$. We have thus the equality of $2$-periodic distributions
\begin{equation}\label{TGRH}- \frac{2}{\pi} ( {\rm sin} \pi x ) {\rm T}_H \,=\, 2 \sum_{n \equiv 1 \bmod 2} \frac{e^{i \pi n x}}{i\pi n}.\end{equation}
We recognize on the right-hand side the primitive with sum $0$ of $\delta_0-\delta_1$: this is the $2$-periodic (distribution) function which is $1$ on $]0,1[$ and $-1$ on $]-1,0[$. We conclude by multiplying \eqref{TGRH} with ${\rm sin}\,\pi x$.\end{pf}

By this lemma, the inner product on $\mathcal{P}^+$ which occurs in Proposition \ref{orthpol} is $ \frac{\pi}{2}\, \int_0^1 \,P(x)\, Q(x)\, {\rm sin}\, \pi x\, {\rm d}x$. Recall that any $P \in \mathcal{P}^+ \cap \mathcal{P}_n$ may be written uniquely $P(x)=\underline{P}({\rm cos} \,\pi x)$ with $\underline{P} \in \R[T]$ of degree $\leq n$ and $\underline{P}(-T) = (-1)^n \underline{P}(T)$. Using the change of variables $u= {\rm cos} \, \pi x$, the scalar product above is nothing else than $\frac{1}{2} \int_{-1}^1 \underline{P}(u) \underline{Q}(u) {\rm d} u$. The orthogonal polynomials ${\rm L}_n$ of this scalar product on $\R[T]$ normalized by ${\rm L}_n(1)=1$ are the famous {\it Legendre polynomials}. By Proposition~\ref{orthpol} and the classical generating series of Legendre polynomials we deduce 
\begin{equation} \label{legendre} \frac{1}{\sqrt{1 - e^{-i\pi x} z}} \frac{1}{\sqrt{1-e^{i\pi x} z}} = \sum_{n \geq 0} {\rm P}_n^H(x) \,z^n, \, \, \,\, \, H(t)=e^{-\frac{t}{2}}.\end{equation}

\noindent Denote by ${\rm h}_n=1+\frac{1}{2}+\dots+\frac{1}{n}$ the harmonic series, and set ${\rm h}_0=0$. 

\begin{prop}\label{calctnGRH} Assume $H(t)=e^{-\frac{t}{2}}$. Then $\underline{{\rm P}_n^H}$ is the $n$-th Legendre polynomial and we have ${\rm t}_n^H = \psi(\frac{1}{2})+ {\rm h}_n$ for all $n\geq 0$.
\end{prop}

\begin{pf} It only remains to prove the assertion about ${\rm t}_n^H$. We have {\small
\begin{equation} \,\label{explpnh0} e^{i \pi n x}\,{\rm P}_n^H(x) \,= \,\sum_{k=0}^n 4^{-n} {{2k}\choose{k}} {{2(n-k)}\choose{n-k}}e^{2i\pi k x} \end{equation}}
\par\noindent by \eqref{legendre}, the {\it binomial identity} $(1-z)^\alpha = \sum_{n \geq 0} (-1)^n{{\alpha}\choose{n}} z^n$ for $\alpha =- \frac{1}{2}$, and the equality $(-1)^k{{-\frac{1}{2}}\choose{k}}=4^{-k} {{2k}\choose{k}}$ for $k\geq 0$ an integer. Using the formula ${\rm t}_n^H = \langle {\rm T}_H,e^{i \pi n x}{\rm P}_n  \rangle =  \langle {\rm T}_H,e^{i \pi n x}{\rm P}_n - {\rm P}_n(0) \rangle + {\rm P}_n(0) \langle {\rm T}_H, 1  \rangle $ given by Proposition \ref{orthpol}, the identity ${\rm P}_n(0)=1$, as well as the tautology $\langle {\rm T}_H, e^{2i\pi kx}  \rangle = \psi(\frac{1}{2}+|k|)$, the identity to be proved is reduced to the following equality of formal power series of $z$:
{\small $$-\,(1-z)^{-1}{\rm log}(1-z)\,= \,(1-z)^{-\frac{1}{2}}\sum_{k\geq 0} (-1)^k{{-\frac{1}{2}}\choose{k}} (\,\psi(\frac{1}{2}+k)-\psi(\frac{1}{2})\,) z^k.$$}

\noindent But if we multiply this identity by $-(1-z)^{\frac{1}{2}}$ on both sides, we recognize the derivative at $\alpha=-\frac{1}{2}$ of the binomial identity, so we are done. \end{pf}
\noindent Recall that we have $\psi(\frac{1}{2})\,=\,-\gamma\,-\,{\rm log} \,4$ where $\gamma=-\Gamma'(1)$ is the {\it Euler-Mascheroni constant}. 
\begin{cor}\label{corthgrh} Assume $H(t)=e^{-\frac{t}{2}}$. The quadratic form $({\rm log}\, 2\pi) \, \phi_n^2 - {\rm q}_n^H$ is positive definite if and only if we have $0 \leq n \leq 24$. For $n \geq 25$, this form is negative on the vector $({{2k}\choose{k}} {{2(n-k)}\choose{n-k}})_{0 \leq k \leq n}$. \end{cor}

\begin{pf} By Proposition \ref{crittnh} (i), the quadratic form $({\rm log}\, 2\pi) \, \phi_n^2 - {\rm q}_n^H$ is positive definite if and only if we have ${\rm log} \, 2 \pi > {\rm t}_n^H$. By Proposition \ref{calctnGRH} the inequality ${\rm log} \, 2 \pi > {\rm t}_n^H$ is equivalent to 
$${\rm h}_n < \gamma + {\rm log} \, 8\pi.$$
But this holds if and only if $n \leq 24$: up to $10^{-2}$, the right-hand side is $\approx 3.80$ and we have ${\rm h}_{24} \approx 3.78$ and $\frac{1}{25}=0.04$. This also shows ${\rm t}_n^H > {\rm log}\, 2\pi$ for $n\geq 25$, and concludes the proof by Formula \ref{explpnh0}. \end{pf}

\subsection{The case $H(t)\,=\,e^{-\frac{t}{2}}\,{\rm cosh}(\frac{t}{2})^{-1} \,=\, 2\, e^{-t} \,(1+e^{-t})^{-1}$}\label{casHsansgrh} In this case we easily see $\psi_H(s)\,=\, {\rm log}\, 2\, + \psi(\frac{1+s}{2})$. For any integer $m$, we get 
\begin{equation} \label{ratlog2psih}  \psi_H(|m|) \,+\, \gamma + \,(-1)^{m}\, {\rm log}\, 2\, \, \, \in \Q, \end{equation} thanks to the classical formulae $\psi(s+1)\,=\,\frac{1}{s}+\psi(s)$, $\psi(1)\,=\,-\,\gamma$ and $\psi(\frac{1}{2})\,=\,-\,\gamma \,-\, {\rm log} \,4.$ In particular, by Formula (iii) of Proposition \ref{defth} applied to $t\,=\,-\,\gamma$, the element $1/(\gamma+{\rm t}_n^H)$ is the sum of all the coefficients of the inverse of the matrix $(\gamma+\psi_H(|i-j|))_{0 \leq i,j \leq n}$, hence an element of $\Q(u)$ with $u={\rm log}\, 2$. On our personal computer, \texttt{PARI} computes this element in a few seconds for $n\leq 25$; we find for instance: 
\vspace{-.4cm}
\begin{table}[h!]\renewcommand{\arraystretch}{1.5}
{\begin{tabular}{c||c|c|c|c|c|c|c|c}
 $n$  & $0$ & $1$ & $2$ & $3$ & $4$ & $5$ & $6$ & $7$ \\ 
\hline 
${\rm t}_n^H+\gamma$ & $-u$ & $0$ & $\frac{u}{4u-1}$ & $\frac{2}{3}$ & $\frac{107 u - 48}{128u - 59}$ & $\frac{153}{145}$ & $\frac{52333 u - 27852}{44292 u - 23701}$ & $\frac{446591}{335349}$ \\
\end{tabular}
\label{tabletnogrh}
}
\end{table}

\noindent Note that $\gamma+{\rm t}_n^H$ is always the quotient of two elements of $\Q \, + \, \Q\, u$. Indeed, both the determinant and the cofactors of $(\gamma+\psi_H(|i-j|))_{0 \leq i,j \leq n}$ are elements of $\Q + \Q\, u$, by Formula \eqref{ratlog2psih} and the fact that the matrix $((-1)^{i+j})_{0 \leq i,j \leq n}$ has rank $1$. Moreover, the vector ${\rm v}_n^H$ is a generator of the kernel of the matrix $(\psi_H(|i-j|)+\gamma-(\gamma+{\rm t}_n^H))_{0 \leq i,j \leq n}$ by Proposition \ref{crittnh} (i), hence it has coordinates in $\Q(u)$ (and we can say more: see Corollary \ref{calctnSANSGRH}). \texttt{PARI} computes again this element in a few seconds for $n\leq 25$. This computation shows in particular that the coefficients of ${\rm v}_n^H$ are positive for $n\leq 25$; we expect, but did not prove, that this holds for all $n$. \ps\ps

\begin{cor} \label{corthsansgrh} Assume $H(t)=e^{-\frac{t}{2}}{\rm cosh}(\frac{t}{2})^{-1}$. The quadratic form $({\rm log}\, 2\pi) \, \phi_n^2 - {\rm q}_n^H$ is positive definite if and only if we have $0 \leq n \leq 23$. For $n=24,25$, this form is negative on the vector $x \in \R^{n+1}$ defined by $x_0=x_{n}=4$, $x_k=1$ for $0<k<n$ and $k \neq \frac{n}{2}$, and $x_k=0$ for $k=\frac{n}{2}$. \ps

\end{cor}

\begin{pf}  The quadratic form of the statement is positive definite if and only if we have ${\rm log} \, 2 \pi > {\rm t}_n^H$ by Proposition \ref{crittnh} (i). Up to $10^{-3}$, we have ${\rm log} \, 2 \pi \approx 1.838$. Plugging $u={\rm log} \, 2$ in the formula for ${\rm t}_n+\gamma$ found above, the computer tell us ${\rm t}_{23}^{H} \approx 1.821$ and ${\rm t}_{24}^{H} \approx 1.861$. This proves the first assertion of the corollary as ${\rm t}_n^{H}$ is a nondecreasing function of $n$ by Proposition \ref{defth} (iii). For the second assertion, the computer tells us that for $y=\frac{x}{\phi_n(x)}$, we have ${\rm q}_n^H(y) \approx 1.852$ for $n=24$ and ${\rm q}_n^H(y) \approx 1.885$ for $n=25$, up to $10^{-3}$ .
\end{pf}
\ps\ps
{\footnotesize
We are not aware of any simple close formula for ${\rm t}_n^H$, nor for the vector ${\rm v}_n^H$. Let us briefly discuss what Proposition \ref{orthpol} tells us. Denote by $\{1-2x\}$ the $1$-periodic (distribution) function which coincides with $1-2x$ on $]0,1[$. An argument similar to the one of Lemma \ref{sinTHgrh} shows
$$ - \,{\rm sin} \,2 \pi x\, \, \,{\rm T}_H \,=\, 2\pi  \sum_{n \in \Z-\{0\}} \frac{e^{2i\pi nx}}{2 i\pi n} \,=\,\pi \{1-2x\}.$$
As ${\rm T}_H$ is nonsingular at $\frac{1}{2}$ by Lemma \ref{identTH}, this implies the equality of $2$-periodic distributions $- {\rm sin} \pi x\, \, {\rm T}_H \,= \,\frac{\pi}{2} \frac{\{1-2x\}}{{\rm cos}\, \pi x}$, hence the following equalities:
{\scriptsize $$ - \langle ({\rm sin} \pi x)^2\, {\rm T}_H, PQ \rangle \,=\, \frac{\pi}{2} \int_0^1 \,P(x)\,Q(x) \,(1-2x)\,{\rm tan} \,\pi x \,{\rm d}x \,= \,\int_{-1}^1 \,\underline{P}(v) \,\underline{Q}(v) \,\frac{{\rm arcsin} \,v}{\pi \, v} \,{\rm d} \,v$$}
\par \noindent for any $P,Q \in \mathcal{P}^+$. We have not found this {\it weight} ${\rm w}(v)=\frac{{\rm arcsin} \,v}{\pi\, v}$ studied in the literature of orthogonal polynomials. Its {\it moments} $\mu_k=\int_{-1}^1 v^k {\rm w}(v){\rm d}v$ obviously vanish for $k\geq 1$ odd. Using ${\rm arcsin}'(v) = (1-v^2)^{-\frac{1}{2}}$, an integration by parts, and Euler's {\it beta} function ${\rm B}(x,y)$, we leave as an exercise to the reader to check the equality $\mu_0= -\frac{1}{2\pi}\,{\rm B}(\frac{1}{2},\frac{1}{2})(\psi(\frac{1}{2})-\psi(1))={\rm log}\, 2$, and the formula $\mu_k = \frac{1}{\pi k}(\pi-{\rm B}(\frac{k+1}{2},\frac{1}{2}))=\frac{1}{k}( 1-  2^{-k} {{k}\choose{k/2}})$ and for $k>0$ even (a rational number). As we have ${\rm w}(-v)=-{\rm w}(v)$, we obtain a family of orthogonal polynomials $Q_0,Q_1,\dots,Q_m$ for ${\rm w}$ using the following obvious formulas (compare with \cite[Chap. II \S 2.2\,\,(2.2.6)]{szego}), 
{\scriptsize $$Q_{2n}(v) = \left| \begin{array}{ccccc} 
\mu_0 & \mu_2 & \cdots & \mu_{2n} \\ 
\mu_2 & \mu_4  & \cdots & \mu_{2n+2} \\
\cdots & \cdots & \cdots & \cdots \\
\mu_{2n-2} & \mu_{2n} &  \cdots & \mu_{4n-2} \\
1 & {\it v}^2 & \cdots & {\it v}^{2n} \\
\end{array} \right|, \, \, \, \, \, 
Q_{2n+1}(v) = \left| 
\begin{array}{ccccc} 
\mu_2 & \mu_4 & \cdots & \mu_{2n+2} \\ 
\mu_4 & \mu_6 & \cdots & \mu_{2n+4} \\
\cdots & \cdots &  \cdots & \cdots \\
\mu_{2n} & \mu_{2n+2} & \cdots & \mu_{4n} \\
{\it v} & {\it v}^3 & \cdots & {\it v}^{2n+1} \\
\end{array} \right|.
$$  }

\par \noindent It is clear on these formulas that the coefficients of $Q_m$ are in $\Q$ for $m$ odd, and in $\Q\,+\,\Q \,{\rm log}\, 2$ for $m$ even (with $Q_m(0) \in \Q$). A similar property holds thus for the number $Q_m(1)$, and by Proposition \ref{orthpol} we have $Q_m(1) \neq 0$ and $\underline{{\rm P}_m^H}\,=\, Q_m/Q_m(1)$. \ps

\begin{cor}\label{calctnSANSGRH} For $n$ odd, we have ${\rm t}_n^H+\gamma \, \in \, \Q$ and ${\rm v}_n^H \,\in\, \Q^{n+1}$. For $n$ even, the coefficients of ${\rm v}_n^H$ are the quotient  two elements of $\Q \, + \, \Q\, {\rm log}\,2$.\end{cor}

\begin{pf} We have just proved the assertions about ${\rm v}_n^H$. Assume $n$ odd. We have ${\rm t}_n^H = \langle {\rm T}_H, e^{i \pi n x} {\rm P}_n^H \rangle$ by Proposition \ref{orthpol}. For $0 \leq k \leq n$ odd, we have $\frac{n+k}{2} \not \equiv \frac{n-k}{2} \bmod 2$, so the element $\langle {\rm T}_H, e^{i \pi nx}\, {\rm cos} \,k\pi x \rangle\,=\, \frac{1}{2}(2 \,{\rm log} \,2\, + \psi(\frac{1}{2}+\frac{n+k}{4})+\psi(\frac{1}{2}+\frac{n-k}{4}))$ is in $-\gamma  + \Q$. We conclude as ${\rm P}_n^H$ as rational coefficients and satisfies ${\rm P}_n^H(0)=1$.  \end{pf}

}
\section{Archimedean preliminaries}\label{arch} \subsection{Grothendieck rings of Archimedean Weil groups} \label{archgen} We fix an Archimedean local field $E$. There is a canonical embedding $\R \rightarrow E$, and we either have $E=\R$, or $E \simeq \C$ and ${\rm Aut}(E)=\Z/2\Z$. The {\it Weil group of $E$} is a topological group ${\rm W}_E$ that may be defined as follows \cite{tate}: we set ${\rm W}_E = E^\times$ if $E$ is complex, and ${\rm W}_\R = \C^\times  \,\coprod \,j \C^\times$, where $j^2$ is the element $-1$ of $\C^\times$ and with $jzj^{-1}=\overline{z}$ for all $z \in \C^\times$. If $E$ is complex, the choice of an isomorphism $E \isomo \C$ allows to identify ${\rm W}_E$ with ${\rm W}_\C=\C^\times$.\ps\ps

We recall some basic facts about the complex linear representations of ${\rm W}_E$ (always assumed finite dimensional and continuous). First of all, there is a unique homomorphism $|\cdot|_E : {\rm W}_E \rightarrow \R_{>0}$ sending any $z \in \C^\times$ to $z\overline{z}$. The natural subgroup $\R_{>0}$ of ${\rm W}_E$ is central, closed, with compact quotient: that quotient is a circle $\mathbb{S}^1$ when $E$ is complex, an extension of $\Z/2\Z$ by $\mathbb{S}^1$ when $E$ is real. If $V$ is an irreducible representation of ${\rm W}_E$, there is a unique $s\in \C$ such that $V \otimes |\cdot|_E^s$ is trivial on $\R_{>0}$. We denote by ${\rm K}({\rm W}_E)$ the Grothendieck ring of ${\rm W}_E$ and by ${\rm K}_E \subset {\rm K}({\rm W}_E)$ the subring of virtual representations which are trivial on $\R_{>0}$. \ps\ps

We obviously have ${\rm K}_\C = \Z[ \eta, \eta^{-1} ]$ where $\eta : \C^\ast \rightarrow \C^\ast$ is the unitary character defined by $\eta(z)=z/|z|_\C^{1/2}$. Moreover, up to isomorphism the irreducible representations of ${\rm W}_\R$ which are trivial on $\R_{>0}$ are 
$$1, \epsilon_{\C/\R}, \, \, \,\,{\rm and}\, \, \,\,{\rm I}_w \,\,{\rm for} \,\,w>0,$$
where $1$ is the trivial character, $\epsilon_{\C/\R}$ is the order $2$ character, and with ${\rm I}_w \,=\, {\rm Ind}_{\C^\times}^{{\rm W}_\R}\, \eta^w$ for $w \in \Z$. We also have ${\rm I}_0 \simeq 1 \oplus \epsilon$ and ${\rm I}_w \simeq {\rm I}_{-w}$. The multiplicative structure on ${\rm K}_\R$ follows from the identities 
\begin{equation} \label{ringkr} {\rm I}_w \otimes {\rm I}_{w'} \simeq {\rm I}_{w+w'} \oplus {\rm I}_{w-w'} \, \, \, \, {\rm and}\, \, \, \, {\rm I}_w \otimes \epsilon_{\C/\R} \simeq {\rm I}_w.\end{equation}

\begin{definition} Let $w \geq 0$ be an integer. We define ${\rm K}_\C^{\leq w} \subset {\rm K}_\C$ as the subgroup generated by the $\eta^v$ with $|v| \leq w$ and $v \equiv w \bmod 2$. We also define ${\rm K}_\R^{\leq w}$ as the subgroup of ${\rm K}_\R$ generated by the ${\rm I}_v$ for $0 \leq v\leq w$ and $v \equiv w \bmod 2$, and also by $1$ in the case $w$ is even; this is a free abelian group of rank $\frac{w+1}{2}$ when $w$ is  odd, $\frac{w}{2}+2$ when $w$ is even.
\end{definition}

Duality on representations induces a ring involution $U \mapsto U^\vee$ of ${\rm K}_E$. This involution sends $\eta$ on $\eta^{-1}$, and is trivial on ${\rm K}_\R$. 
Moreover, induction and restriction define two natural $\Z$-linear maps ${\rm ind} : {\rm K}_\C \rightarrow {\rm K}_\R$ and ${\rm res} \, \, {\rm K}_\R \rightarrow {\rm K}_\C$, which both preserve the ``$\leq w$ subgroups''. For $U \in {\rm K}_\C$ and $V \in {\rm K}_\R$ we have the trivial formulas
\begin{equation}\label{indres} {\rm res} \, \, {\rm ind}\, U \, \, =\, \, U+U^\vee\, \, \, {\rm and} \, \, \, {\rm ind}\,\, {\rm res}\, V  \,\,=\,\, {\rm I}_0\,\cdot \, V.\end{equation}
\ps

\subsection{Some quadratic forms on ${\rm K}_E$.}\label{quadgamma} Let $E$ be as in \S \ref{archgen} and $\Gamma(s)$ be the classical $\Gamma$ function. For any $U \in {\rm K}({\rm W}_E)$, there is meromorphic function $s \mapsto \Gamma(s,U)$, defined according to Tate's recipe \cite{tate} as the ``$\Gamma$ factor'' of $U \otimes |.|_E^s$. It satisfies $\Gamma(s,V \oplus W)=\Gamma(s,V)\Gamma(s,W)$ for all $V,W \in {\rm K}({\rm W}_E)$, and $\Gamma(s,U) = \Gamma(s,{\rm ind}\, U)$ for all $U \in {\rm K}_\C$. It is thus enough to give $\Gamma(s,U)$ for $U$ in a $\Z$-basis of ${\rm K}_\R$, namely:
{\small $$\Gamma(s,1)=\pi^{-\frac{s}{2}}\Gamma(\frac{s}{2})\, \,\, \, {\rm and} \, \, \, \, \, \Gamma(s,{\rm I}_{w}) = 2 (2\pi)^{-s-\frac{w}{2}}\Gamma(s+\frac{w}{2})\, \, \, {\rm for\,\, all}\,\,\, w \geq 0.$$}\par

\noindent Let us fix an integrable function $F : \R \rightarrow \R$ such that the function $t \mapsto \frac{F(t)-F(0)}{t}$ is of bounded variation on $\R-\{0\}$. Under these assumptions, \cite[Prop. 3]{poitou2} shows that for all $U \in {\rm K}_E$ the integral\footnote{Our convention for the Fourier transform of a function $F \in {\rm L}^1(\R)$ is $\widehat{F}(\xi)\,=\,\int_\R F(x) e^{-2i\pi x \xi}\, {\rm d}x$.} 
\begin{equation}\label{defJFEU}{\rm J}_F^E(U) =  -  \int_\R \,\,\frac{\Gamma'}{\Gamma}(\,\frac{1}{2}\,+\,2\pi t \,i,U) \,\widehat{F}(t)\, {\rm d}t\end{equation}
converges, as well as the following concrete formulas:
{\small
\begin{equation} \label{formulesJF} \left\{  \begin{array}{l} {\rm J}_F^\C(\eta^w)={\rm J}_F^\R({\rm I}_w)= F(0) \, {\rm log}\, 2\pi\, -\psi_F(\frac{1+|w|}{2}) \,\, \, \, {\rm for}\, \, {\rm all}\, \, w \in \Z,\\ \\ {\rm J}_F^\R(1-\epsilon_{\C/\R}) \,= \,\int_0^\infty \,F(t)\frac{e^{-\frac{t}{2}}}{1+e^{-t}}{\rm d}t. \end{array} \right.\end{equation}
}
\par \noindent As we will recall later, the integrals \eqref{defJFEU} typically show up in the ``explicit formulas'', which explains their study here. Note that ${\rm J}_F^E$ is a linear map ${\rm K}_E \rightarrow \R$ with the following properties:  
\begin{equation} \label{propJFE} {\rm J}_F^\R \circ {\rm ind} = {\rm J}_F^\C,  \hspace{.5 cm}{\rm  hence}\hspace{.5 cm} {\rm J}_F^E(U^\vee)={\rm J}_F^E(U)\,\, {\rm for}\,\,\,{\rm all}\,\,\,U \in {\rm K}_E.\end{equation} 
\noindent Observe also that the formulas \eqref{formulesJF} make sense more generally for an arbitrary function $F : [0,+\infty[ \rightarrow \R$ such that $H(t) := F(t)e^{-\frac{t}{2}}$ satisfies the assumptions of \S \ref{somequad}, if we interpret $\psi_F(\frac{1+|w|}{2})$ as $\psi_H(\frac{|w|}{2})$. We fix now such a function $F$, and denote by ${\rm J}_F^E : {\rm K}_E \rightarrow \R$  the unique collection of linear maps satisfying \eqref{formulesJF} and \eqref{propJFE}. For every $E$, we define a symmetric bilinear form on ${\rm K}_E$ by the formula
\begin{equation} \label{defUVFE} \langle U, V \rangle_F^E = {\rm J}_F^E(U \cdot V^\vee).\end{equation}
We also denote by $\dim : {\rm K}_E \rightarrow \Z$ the ring homomorphism giving the virtual dimension. The following immediate lemma is the bridge between sections \ref{sect1} and \ref{arch}. \ps\ps

\begin{lemme} \label{matcompl} Fix a function $F : [0,+\infty[ \rightarrow \R$ and assume that $H(t)=F(t)e^{-t/2}$ satisfies the assumptions of \S \ref{somequad}. For any $c \in \R$ and any integer $w \geq 0$, the Gram matrix of the bilinear form $\langle U, V \rangle_F^\C \,-\, c\, \dim U\, \dim V$ in the $\Z$-basis $\eta^{-w}, \eta^{-w+2}, \dots, \eta^{w-2}, \eta^{w}$ of ${\rm K}_\C^{\leq w}$ is $$( H(0){\rm log}\, 2\pi\, - c  - \psi_H(|i-j|))_{0 \leq i,j \leq w}.$$ 
\end{lemme}
\ps\ps
\noindent The Gram matrix of $\langle U, V \rangle_F^\R$ in the natural $\Z$-basis of ${\rm K}_\R^{\leq w}$ is not as nice as the one above, because of the relations \eqref{ringkr}. Nevertheless, the following proposition reduces the study of the {\it real} case to that of the {\it complex} one.\ps\ps

\begin{prop} \label{wrwc} Fix a function $F : [0,+\infty[ \rightarrow \R_{\geq 0}$ and assume that $H(t)=F(t)e^{-t/2}$ satisfies the assumptions of \S \ref{somequad}, with $H$ nonzero in ${\rm L}^1(]0,+\infty[)$. For any $t \in \R$ and any integer $w \geq 0$, there is an equivalence between: \ps\ps
\noindent (i) the form $U \mapsto \langle U, U \rangle_F^\R - \,\frac{t}{2}\, ({\dim U})^2$ is positive definite on ${\rm K}_\R^{\leq w}$,\ps \ps
\noindent (ii) the form $U \mapsto \langle U, U \rangle_F^\C \,-\, t \,({\dim U})^2$ is positive definite on ${\rm K}_\C^{\leq w}$.\ps
\end{prop} 

\begin{pf} Fix $w$ and $t$ as in the statement and consider the $\Q$-linear map ${\rm K}_\R^{\leq w} \otimes \Q \rightarrow {\rm K}_\C^{\leq w} \otimes \Q$ induced by ${\rm res}$. Its kernel $K$ is generated by $1-\epsilon_{\C/\R}$ for $w$ even, and vanishes for $w$ odd. 
Moreover, we have ${\rm K}_\R^{\leq w} \otimes \Q = K \oplus I$, where $I$ is the $\Q$-linear span of the ${\rm I}_v$ for $0 \leq v \leq w$ and $v \equiv w \bmod 2$. The image of ${\rm res}$ is the subspace of elements $U$ in ${\rm K}_\C^{\leq w} \otimes \Q$ with $U^\vee=U$.

\ps
We first claim that the quadratic form of assertion (i) is positive definite on $K$, and that it is positive definite on the whole of ${\rm K}_\R^{\leq w}$ if, and only, it is so on the subspace $I$. Indeed, note that $K$ and $I$ are orthogonal with respect to both bilinear forms $\langle U, V \rangle_F^\R$ and $\dim \,\,U \cdot V$, by the identity $KI=0$ in  ${\rm K}_\R$ (an immediate consequence of Formulas \eqref{ringkr}). Moreover, under the assumption on $F$, Formula \eqref{formulesJF} shows that we have ${\rm J}_F(1-\epsilon_{\C/\R}) >0$. As $U \mapsto (\dim U)^2$ obviously vanishes on $K$, this proves the claim. \ps
Next, we claim that the linear map ${\rm res}$ defines an isometric embedding $I \rightarrow {\rm K}_\C^{\leq w}$, where the source is equipped with twice the form of assertion (i), and the target with that of assertion (ii). Indeed, we have $\dim \,\circ\, {\rm res} \,= \,\dim$ and $$\langle\, {\rm res}\, U,\,{\rm res} \,V \,\rangle_F^\C \, =\,{\rm J}_F^\R({\rm ind}\, {\rm res}\, \,U \cdot V^\vee)\, = \,2 \,\langle U, V \rangle_F^\R$$ for all $U,V \in I$. The first equality above is a consequence of \eqref{propJFE}, and the second one follows from the formula ${\rm ind}\, {\rm res} \,W\, =\, {\rm I}_0\, \cdot W\, =\, 2\, W$, which holds by \eqref{ringkr} and \eqref{indres} for all $W$ in ${\rm K}_\R^{\leq w}$ with $W = W \cdot\epsilon_{\C/\R}$, applied to the element $W=U \cdot V^\vee$. In particular, we have proved that assertion (ii) implies assertion (i). \ps
Last, under our assumption on $F$, Lemma \ref{matcompl} and Proposition \ref{crittnh} show that the form of assertion (ii) is positive definite on the hyperplane $\dim =0$ of ${\rm K}_\C^{\leq w}$. Better, by Proposition \ref{crittnh} (i) we know that this form is positive definite on the whole of ${\rm K}_\C^{\leq w}$ if, and only if, it is $>0$ on the unique element $V \in {\rm K}_\C^{\leq w} \otimes \R$ with $\dim V=1$ which is orthogonal (with respect to $\langle -, - \rangle_F^\C$)  to the hyperplane $\dim =0$. By uniqueness, this element verifies $V=V^\vee$, hence lies in the image of $I \otimes \R$ under the isometric embedding ${\rm res} \otimes \R$. This shows that assertion (i) implies assertion (ii), and concludes the proof of the proposition.
\end{pf}

\subsection{Algebraic Harish-Chandra modules} \label{archalg} We still assume $E$ is an Archimedean local field. The local Langlands correspondence for ${\rm GL}_n(E)$ is a natural bijection $U \mapsto {\rm L}(U)$ between the set of isomorphism classes of irreducible admissible Harish-Chandra modules for ${\rm GL}_n(E)$ (say, with respect to the standard maximal compact subgroup) and the set of isomorphism classes of $n$-dimensional semisimple representations of ${\rm W}_E$, for all integers $n\geq 1$ \cite{knappmotives}.  For any irreducible Harish-Chandra module $U$ of ${\rm GL}_n(E)$, we may and shall view ${\rm L}(U)$ as an (effective) element of ${\rm K}({\rm W}_E)$. \ps\ps

Let $U$ be an irreducible Harish-Chandra module for ${\rm GL}_n(E)$. Recall that $U$ has an infinitesimal character, that we may view following Harish-Chandra and Langlands as a semisimple conjugacy class in ${\rm M}_n(E \otimes_\R \C)$. This is really the datum of one or two semisimple conjugacy classes ${\rm inf}_\sigma\, U$ in ${\rm M}_n(\C)$, indexed by the embeddings $\sigma : E \rightarrow \C$. The eigenvalues of ${\rm inf}_\sigma\, U$ will be called the {\it weights of $U$ relative to $\sigma$}. The module $U$ is said {\it algebraic} if all of its weights are in $\Z$ (relative to each $\sigma$). This property may be read off from ${\rm L}(U)$ as follows. Fix an embedding $\sigma : E \rightarrow \C$; it gives rise to a natural inclusion $\C^\times \rightarrow {\rm W}_E$. The module $U$ is algebraic if, and only if, there are integers $p_{i,\sigma},q_{i,\sigma} \in \Z$ for $i=1,\dots,n$ with
\begin{equation} \label{relLweight} {\rm L}(U)_{|\C^\times}  \simeq \bigoplus_{i=1}^n\,\,(z \mapsto  z^{p_{i,\sigma}} \overline{z}^{q_{i,\sigma}}).\end{equation}
This property does not depend on the choice of $\sigma$; moreover, the multiset of the $p_{i,\sigma}$ (resp. $q_{i,\sigma}$) above coincides with the multiset of weights of $U$ relative to $\sigma$ (resp. $\overline{\sigma}$). Note in particular that those two multisets are clearly the same for $E=\R$, by the relation $j z j^{-1} = \overline{z}$ in ${\rm W}_\R$. \ps

\subsection{Algebraic automorphic representations} \label{archalgglob} Assume now that $E$ is a number field and let $\pi$ be a cuspidal automorphic representation of ${\rm GL}_n$ over $E$. We say that $\pi$ is algebraic if so is $\pi_v$ for each $v \in {\rm S}_\infty(E)$, where ${\rm S}_\infty(E)$ denotes the set of Archimedean places of $E$. Assume $\pi$ is algebraic. For each field embedding $\sigma : E \rightarrow \C$, with underlying place $v \in {\rm S}_\infty(E)$, we have by \eqref{relLweight} a collection of weights $(p_{i,\sigma},q_{i,\sigma})$ associated with $\pi_v$ and $\sigma$. We will say that $\pi$ is {\it effective} if we have $p_{i,\sigma}, q_{i,\sigma} \geq 0$ for all $i$ and $\sigma$. As observed by Clozel  \cite[Lemma 4.9]{clozel}, the Jacquet-Shalika estimates together with the Dirichlet units theorem impose strong constraints on the couples $(p_{i,\sigma},q_{i,\sigma})$, namely that the quantity $w=p_{i,\sigma}+q_{i,\sigma}$ depends neither on $i$, nor on $v$ or $\sigma$. This $w$ will be called the {\it motivic weight}\footnote{Some authors (including us) sometimes rather call $-w$ the motivic weight of $\pi$.} of $\pi$ and denoted by ${\rm w}(\pi)$. 
For each $\sigma : E \rightarrow \C$, with underlying place $v \in {\rm S}_\infty(E)$, we have by \eqref{relLweight} an isomorphism $ {\rm L}(\pi_v \otimes |\det|_v^{-{\rm w}(\pi)/2})_{|\C^\times} \, \, \simeq \, \, \bigoplus_{i=1}^n \, \, \eta^{2\, p_{i,\sigma}-{\rm w}(\pi)}$, hence:  \ps\ps

\begin{lemme}\label{lemmemult} Let $\pi$ be a cuspidal automorphic representation of ${\rm GL}_n$ over the number field $E$. Assume $\pi$ is algebraic, effective, of motivic weight $w$. Then for any $v \in {\rm S}_\infty(E)$ the element ${\rm L}(\pi_v \otimes |\det|_v^{-w/2})$ lies in ${\rm K}_{E_v}^{\leq w}$. 
\end{lemme}

\noindent Note that twisting a cuspidal algebraic $\pi$ globally by $|\det|^m$ with $m \in \Z$ preserves algebraicity, translates all the weights by $m$, and the motivic weight by $2m$.



\section{The explicit formula for Rankin-Selberg ${\rm L}$-functions} \label{sectrs}

\subsection{Rankin-Selberg $L$-functions} \label{RS}  Let $E$ be a number field and $\pi$ a cuspidal automorphic representation of ${\rm GL}_n$ over $E$. Consider the associated Rankin-Selberg ${\rm L}$-function ${\rm L}(s,\pi \times \pi^\vee)$ defined by Jacquet, Piatetski-Shapiro and Shalika  \cite{js,jps1}. The analytic properties of this ${\rm L}$-function, such as the absolute convergence of its Euler product for ${\rm Re} \,s >1$, its meromorphic continuation to $\mathbb{C}$ with single and simple poles at $0$ and $1$, and its functional equation $s \rightarrow 1-s$, proved by those authors (with a contribution of Moeglin and Waldspurger for the assertion about the poles \cite{moeglinwaldspurger}), will be used here in a crucial way, as well as the boundedness in vertical strips away from the poles proved by Gelbart and Shahidi in \cite{gs}. See \cite{cogdell} for a survey of these results. Set
$$\Lambda(s,\pi \times \pi^\vee) \,=\, \left( \prod_{v \in {\rm S}_\infty(E)} \Gamma(s, {\rm L}(\pi_v) \otimes {\rm L}(\pi_v)^\vee)\,\,\right)\, {\rm L}(s,\pi \times \pi^\vee),$$
where $\Gamma(s,-)$ is the $\Gamma$-factor recalled in \S \ref{archgen}. The functional equation reads $$\Lambda(s,\pi \times \pi^\vee ) = \epsilon(1/2,\pi \times \pi^\vee)\,{\rm N}(\pi \times \pi^\vee)^{\frac{1}{2}-s}\,  \Lambda(1-s,\pi \times \pi^\vee)$$ where $\epsilon(\frac{1}{2},\pi \times \pi^\vee) \in \C^\times$ is the global root number and ${\rm N}(\pi \times \pi^\vee)$ is an integer that we shall refer to as the {\it absolute conductor} of the pair $\{\pi, \pi^\vee\}$. This latter number has the form 
\begin{equation} \label{formcondabs} {\rm N}(\pi \times \pi^\vee) \, = \, ({\rm disc} \,E)^{n^2} \, \, ||\, \mathcal{N}(\pi \times \pi^\vee)\,||\end{equation}
where ${\rm disc} \,E$ is the discriminant of $E$, and $\mathcal{N}(\pi \times \pi')$ the {\it Artin conductor} of the pair $\{\pi,\pi^\vee\}$ (a nonzero ideal of the ring $\mathcal{O}_E$ of integers of $E$), and where $||I||$ denotes the (finite) cardinality of $\mathcal{O}_E/I$ for any nonzero ideal $I$ of $\mathcal{O}_E$. \ps\ps
Denote by ${\rm S}_f(E)$ the set of finite places of $E$, and for $v \in {\rm S}_f(E)$ denote by $q_v$ the cardinality of the residue field of $v$. A last important property that we shall use, which is proved in \cite[\S 2, Lemma a]{hram}, is that if we write for ${\rm Re}\, s>1$ $$ - \frac{{\rm L}'}{{\rm L}}(s,\pi \times \pi^\vee)= \sum_{v,k} {\rm a}_{v^k}(\pi \times \pi^\vee)\, \frac{{\rm log}\, q_v}{q_v^{ks}},$$ 
the sum being over all $v \in {\rm S}_f(E)$ and all integers $k\geq 1$, then we have
\begin{equation}\label{posanpipiprime} {\rm a}_{v^k}(\pi \times \pi^\vee) \in \R_{\geq 0}\, \, \, \,  {\rm for}\, \, {\rm  all}\,\, v \in {\rm S}_f(E) \,\,{\rm and}\,\, k\geq 1.\end{equation}

\subsection{The Riemann-Weil explicit formula}\label{testfunctions} The statement of the explicit formula requires the choice of an auxiliary function called a {\it test function}. Following Weil, Poitou \cite{poitou2} and Mestre, we shall mean by this any {\it even function $F : \R \rightarrow \R$} satisfying the following three conditions:\footnote{The reader would not loose much by assuming in this paper that a test function is any even function $F : \R \rightarrow \R$ which is compactly supported and of class $\mathcal{C}^2$.}\par\ps
 {\it 
 {\rm (TFa)} There exists $\epsilon>1/2$ such that the function $x \mapsto F(x)e^{\epsilon |x|}$ is in ${\rm L}^1(\R)$ and of bounded variation. \par 
{\rm (TFb)} $x \mapsto (F(x)-F(0))/x$ is of bounded variation on $\R-\{0\}$.\par
{\rm (TFc)} At each point, $F$ is the arithmetic mean of its left and right limits {\rm (}note that $F$ has bounded variations by {\rm (TFa)} and {\rm (TFb)}{\rm )}.
} \ps

\noindent For any test function $F$ we set $\Phi_F(s) = \widehat{F}(\frac{1-2s}{4i\pi})$, with the convention $\widehat{f}(\xi) = \int_\R f(x) e^{-2i\pi x \xi } dx$ for the Fourier transform of $f$.  Note that $\Phi_F(s)$ is defined and holomorphic in a neighborhood of $0 \,\leq \,{\rm Re} \,s\, \leq 1$. \ps

We fix now a number field $E$ and a cuspidal algebraic automorphic representation\footnote{We reserve the letter $\pi$ here for the area of the unit Euclidean disc.} $\varpi$ of $\GL_n$ over $E$. We assume that $\varpi$ is algebraic of motivic weight $w$ (\S \ref{arch}) and set $$\Lambda(s):=\Lambda(s,\varpi \times \varpi^\vee)=\Lambda(s,(\varpi \otimes |\det|^{-\frac{w}{2}}) \times (\varpi \otimes |\det|^{-\frac{w}{2}})^\vee).$$
The {\it explicit formula associated to $\Lambda(s)$ and to the test function $F$} is the result of the contour integration of $\Phi_F(s) {\rm d} {\rm log} \Lambda(s)$ over the boundary of the rectangle defined by $-\epsilon \leq {\rm  Re}\, s \leq 1+\epsilon$ and $|{\rm Im}\, s|\leq T$ (for suitable $T$'s), with $T \rightarrow + \infty$ and $\epsilon \rightarrow 0$. \ps 

The analytic properties of $\Lambda(s)$ recalled in \S \ref{RS} ensure that Mestre's general formalism \cite[\S 1]{mestre} applies to the pair of functions $\Lambda_1(s)=\Lambda_2(s)={\rm N}(\varpi \times \varpi^\vee)^{\frac{s}{2}} \Lambda(s)$ and ${\rm L}_1(s)={\rm L}_2(s)={\rm L}(s,\varpi \times \varpi^\vee)$ in Mestre's notations.\footnote{Note that contrary to Mestre, we choose to incorporate the natural powers of $\pi$ inside the $\Gamma$ factors. Moreover, $ {\rm L}(s,\pi \times \pi^\vee)$ does not quite satisfies assumption (iv) of Mestre for $c=0$, but the known absolute convergence of this Euler product for ${\rm Re}\, s \, >1$ is all he needs: see \cite[p. 213-214]{mestre} and \cite[p. 2-3]{poitou2}.} Our assumption that $\pi$ is algebraic, which implies ${\rm L}(\varpi_v) \otimes {\rm L}(\varpi_v)^{\vee} \, \in \, {\rm K}_{E_v}$ for each Archimedean place $v$ of $E$ by Lemma \ref{lemmemult}, is only used at the moment to ensure that the $\Gamma$-factors in $\Lambda$ have the form $\Gamma(as+b)$ with $a,b \in \R_{>0}$, hence fit as well Mestre's assumptions. The {\it explicit formula} proved by Mestre in \cite[\S 1.2]{mestre}, a generalization of Weil's formula \cite{poitou2}, is (twice) the equality {\small
$$\sum_{v \in {\rm S}_\infty(E)}\, {\rm J}_F^{E_v}({\rm L}(\varpi_v)\otimes {\rm L}(\varpi_v)^\vee) + \,\,\sum_{v\in {\rm S}_f(E),k\geq 1} {\rm a}_{v^k}(\varpi \times \varpi^\vee)\,F(k {\rm log}\, q_v)\,\,\frac{{\rm log} \,q_v}{q_v^{k/2}}$$ $$=$$ $$\frac{F(0)}{2} \,{\rm log} \, {\rm N}(\varpi \times \varpi^\vee) \,+\, \Phi_F(0)\,-\, \frac{1}{2}\,{\sum_{\{\rho \in \C, \, \Lambda(\rho)=0\}} \,\Phi_F(\rho)\, \, {\rm ord}_{s=\rho}\, \Lambda(s).}$$}
{}\noindent In this formula, ${\rm J}_F^{E_v}$ is the linear map ${\rm K}_{E_v} \rightarrow \R$ introduced in \S \ref{quadgamma} (note that a test function is in ${\rm L}^1(\R)$, so that Formula \eqref{defJFEU} holds) and the sum over the zeros $\rho$ of $\Lambda$ has to be understood as the limit of the finite sums over the $\rho$ with $|{\rm Im} \, \rho | \leq T$ with $T \rightarrow +\infty$. Note that all those zeros verify $0 \,\leq \,{\rm Re}\, \rho \,\leq 1$ by Jacquet and Shalika. \ps

\subsection{The basic inequality}\label{basicinequality}
Let $F$ be a test function. We shall say that $F$ satisfies {\rm (POS)} if the following equivalent properties hold:\ps\ps

(a)  ${\rm Re}\, \widehat{F}(s) \geq 0$ for all $s \in \C$ with $|{\rm Im}\, s |\leq \frac{1}{4\pi}$, \ps

(b)  ${\rm Re}\, \Phi_F(s) \geq 0$ for $0 \leq {\rm Re}\, s \leq 1$, \ps

(c) the function $G(x) = F(x) {\rm cosh} \frac{x}{2}$ satisfies $\widehat{G}(\xi) \geq 0$ for all $\xi \in \R$.\ps\ps

\noindent Indeed, the equivalence of (a) and (b) is obvious, and the one of (b) and (c), observed by Poitou in \cite[p. 6-08]{poitou2}, follows from the trivial equality ${\rm Re}\, \widehat{F}(\pm \frac{i}{4\pi} +\xi) \,=\, \widehat{G}(\xi)$ and the minimum principle. From the explicit formula and \eqref{posanpipiprime}, we obtain the following proposition: 

\begin{prop} \label{inegfex} Let $\varpi$ be a cuspidal algebraic automorphic representation of $\GL_n$ over the number field $E$ and let $F$ be a nonnegative test function. If $F$ satisfies {\rm (POS)} then we have\footnote{\label{footis4p}Note that $\widehat{F}(\frac{i}{4\pi}) = \int_\R F(x) e^{\frac{x}{2}} {\rm d}x$ is a nonnegative real number.}
$$\sum_{v \in {\rm S}_\infty(E)}\, {\rm J}_F^{E_v}({\rm L}(\varpi_v)\otimes {\rm L}(\varpi_v)^\vee)  \leq  \widehat{F}(\frac{i}{4\pi})  +  \frac{1}{2}\,F(0) \, \, {\rm log}\, {\rm N}(\varpi \times \varpi^\vee).$$
If all the zeros of $\Lambda(s,\varpi \times \varpi^\vee)$ satisfy ${\rm Re}\, s = 1/2$, this equality holds if we replace {\rm (POS)} by $\widehat{F}(\xi)\geq 0$ for all $\xi\in \R$.  
\end{prop}

The extension of the Riemann hypothesis that appears in the last assertion of this proposition will appear several times in the sequel, and we shall refer to it as:\ps
{\it {\rm (GRH)} For any cuspidal algebraic automorphic representation $\varpi$ of $\GL_n$ over a number field, the zeros of $\Lambda(s,\varpi \times \varpi^\vee)$ satisfy ${\rm Re}\, s\, = \, \frac{1}{2}$. }
\ps

We will apply Proposition \ref{inegfex} in a classical way to the test functions introduced by Odlyzko. One starts with the real function ${\rm g}$ which is twice the convolution square of the function ${\rm cos}(\pi x)$ on $|x| \leq 1/2$ and $0$ otherwise. This is an even nonnegative function with support in $[-1,1]$, easily seen to be of class $\mathcal{C}^2$, hence a test function. As a consequence, for any real number $\lambda >0$, the functions ${\rm G}_\lambda$ and ${\rm F}_\lambda$ defined by $${\rm G}_\lambda(x)= {\rm g}(x/\lambda) \, \, \, \, {\rm and} \, \, \ \, {\rm F}_\lambda(x) = {\rm G}_\lambda(x){\rm cosh}(x/2)^{-1}$$ are nonnegative test functions as well, with support in $[-\lambda,\lambda]$. By construction we have $\widehat{{\rm g}}\geq 0$ on $\R$, hence $\widehat{{\rm G}_\lambda} \geq 0$ as well, and ${\rm F}_\lambda$ satisfies {\rm (POS)}. We also have ${\rm F}_\lambda(0)={\rm G}_\lambda(0)=1$ and $\widehat{{\rm F}_\lambda}(\frac{i}{4\pi})= \frac{8}{\pi^2} \lambda$.  \ps\ps

\noindent We conclude this paragraph by a lemma that we shall use in \S \ref{chaptopt}. \ps

\begin{lemme}\label{lemmainegF} Let $F$ be a nonnegative test function satisfying {\rm (POS)} {\rm (}resp. ``$\widehat{F}(\xi) \geq 0$ for all $\xi \in \R$''{\rm )}. Then we have $F(x) \leq F(0){\rm cosh}(\frac{x}{2})^{-1}$ {\rm (}resp. $F(x) \leq F(0)${\rm )} for all $x \in \R$. 
\end{lemme}

\begin{pf}  By Poitou's observation recalled above, the function $G(x) = F(x) {\rm cosh}(\frac{x}{2})$ (resp. $G(x)=F(x)$) satisfies $\widehat{G}(\xi) \geq 0$ for all $\xi \in \R$. The (even) function $G$ is summable and of bounded variation on $\R$ by {\rm (TFa)}, so that Fourier inversion applies to it by {\rm (TFc)}: we have $G(x) = \int_\R \widehat{G}(\xi) e^{-2i\pi x \,\xi} {\rm d}\xi$ for $x \in \R$. The inequalities $\widehat{G} \geq 0$ and $G \geq 0$ imply then $G(x) = |G(x)| \leq \int_\R \widehat{G}(\xi) {\rm d}\xi = G(0)$. \end{pf}

\section{Proof of the main theorem}\label{pfmainthm}

\begin{lemme}\label{limfl} For any Archimedean local field $E$, and $U$ in ${\rm K}_E$, we have ${\rm J}_{{\rm F}_\lambda}^E(U) \rightarrow {\rm J}_{{\rm cosh}(\frac{t}{2})^{-1}}^E(U)$ and ${\rm J}_{{\rm G}_\lambda}(U) \rightarrow {\rm J}_{1}(U)$ for $\lambda \rightarrow \infty$.
\end{lemme}
\ps
\begin{pf} As we just explained, the inequality $\widehat{{\rm g}}(\xi) \geq 0$ for $\xi \in \R$ implies by Fourier inversion $0 \leq {\rm g}(x) \leq {\rm g}(0)=1$ for all $x \in \R$. Set $f(x)={\rm cosh}(x/2)^{-1}$. We have thus $0 \leq {\rm G}_\lambda(x) \leq 1$ and $0 \leq {\rm F}_\lambda(x) \leq f(x)$, as well as  ${\rm F}_\lambda(x) \rightarrow f(x)$ and ${\rm G}_\lambda(x) \rightarrow 1$ for $\lambda \rightarrow +\infty$. Formulas \eqref{formulesJF} and a trivial application of the dominated convergence theorem conclude the proof. \end{pf}
\ps
\noindent Remark that neither the function $F=1$, nor $F(t)={\rm cosh}(\frac{t}{2})^{-1}$, is a test function, but for both of them the linear forms ${\rm J}_F^E$ of \S \ref{quadgamma} are well defined, as $H(t)=F(t)e^{-\frac{t}{2}}$ satisfies the assumptions of \S \ref{somequad}.\ps

\begin{definition}\label{deftwrw} Let $w\geq 0$ be an integer. We define ${\rm t}(w)$ {\rm (}resp. ${\rm t}^\ast(w)${\rm )} as the real number ${\rm log} \,2\pi
\, - {\rm t}_w^H$ with $H(t)=e^{-\frac{t}{2}}{\rm cosh}(\frac{t}{2})^{-1}$ {\rm (}resp. $H(t)=e^{-\frac{t}{2}}${\rm )}. We also set ${\rm r}(w)={\rm exp}( {\rm t}(w))$ and ${\rm r}^\ast(w)={\rm exp}( {\rm t}^\ast(w))$. 
\end{definition}

By Proposition \ref{crittnh} (ii), each of those four functions of $w$ is nondecreasing. We also have the inequality ${\rm r}(w) \leq {\rm r}^\ast(w)$ for all $w \geq 0$ by Proposition \ref{inegqh1h2} applied to $H_1(t)=e^{-\frac{t}{2}}{\rm cosh}(\frac{t}{2})^{-1}$ and $H_2(t)=e^{-\frac{t}{2}}$. Moreover, we have the closed formula 
$${\rm r}^\ast(w) \,=\,8\pi\,e^{\gamma-\sum_{1 \leq k \leq w} \frac{1}{k}}\, = \,8 \pi\, e^{-\, \psi(1+w)}\,$$ 
by Proposition \ref{calctnGRH}. The equality between the definition here of ${\rm t}(w)$ and the one of the introduction is a consequence of Proposition \ref{crittnh} (iii). Table \ref{tablerw} in the introduction gives some values of ${\rm r}(w)$ up to $10^{-4}$, we have computed them with \texttt{PARI} using the method described in \S \ref{casHsansgrh}. \ps

Recall that if $E$ is a number field of degree $d$, we denote by ${\rm r}_E = |{\rm disc}\, E|^\frac{1}{d}$ its root-discriminant.\ps
\begin{prop}\label{propmainthm} Let $w\geq 0$ be an integer and $E$ a number field with ${\rm r}_E < {\rm r}(w)$. There are positive definite quadratic forms ${\rm q}_v: {\rm K}_{E_v}^{\leq w} \rightarrow \R$ for $v \in {\rm S}_\infty(E)$, as well as a real number $C>0$, such that for any cuspidal automorphic representation $\pi$ of ${\rm GL}_n$ over $E$ which is algebraic, effective, and of motivic weight $w$, we have 
$$\sum_{v \in {\rm S}_\infty(E)} {\rm q}_{v} ({\rm L}(\pi_v \otimes |\det|_v^{-w/2})) \,\leq C\,+ \,{\rm log} \,\,||\, \mathcal{N}(\pi \times \pi^\vee)\,||.$$
Furthermore, a similar statement also holds with ${\rm r}_E  < {\rm r}^\ast(w)$ if we assume ${\rm (GRH)}$.
\end{prop}

\begin{pf} Let $E$ be a number field and $\varpi$ a cuspidal automorphic representation of ${\rm GL}_n$ over $E$ which is algebraic, effective, of motivic weight $w$.  For each $v \in {\rm S}_\infty(E)$, set $V_v = {\rm L}(\varpi_v \otimes |\det|_v^{-w/2})$; this is an element of ${\rm K}_{E_v}^{\leq w}$ by Lemma \ref{lemmemult}. Of course, we have $\dim V_v=n$ for each $v$, so Formula \eqref{formcondabs} takes the following equivalent form:
{\small $${\rm log}\, {\rm N}(\varpi \times \varpi^\vee)\, =\,{\rm log} \,\,||\, \mathcal{N}(\varpi \times \varpi^\vee) \,||\,\,+\, \sum_{v \in {\rm S}_\infty(E)} (\dim V_v)^2 \, [E_v:\R]\,  \,\,{\rm log}\, {\rm r}_{E}.$$} 
As a consequence, for any nonnegative test function $F$ satisfying {\rm (POS)} (resp. $\widehat{F}\geq0$ under ${\rm (GRH)}$), the basic inequality of Proposition \ref{inegfex} may also be written 
\begin{multline}\label{expleffective} \sum_{v \in {\rm S}_\infty(E)} \left( \langle V_v, V_v \rangle_F^{E_v} -  \, F(0)\,\frac{[E_v:\R]}{2} \, {\rm log}\, {\rm r}_{E}\, \dim V_v^2 \right) \\  \leq  \widehat{F}(\frac{i}{4\pi})  +  \,\frac{F(0)}{2}\, {\rm log}\, \, ||\,\mathcal{N}(\varpi \times \varpi^\vee)\,||.\end{multline}
Assume $F\geq 0$ on $\R$ and $F(0)=1$  (note that $F$ is continuous at $0$), Proposition \ref{wrwc} shows the equivalence between: \ps
(a) $U \mapsto \langle U, U \rangle_F^\R - \,\frac{{\rm log}\, {\rm r}_E}{2} \, {\dim}\, U^2$ is positive definite on ${\rm K}_\R^{\leq w}$,\ps \ps
(b) $U \mapsto \langle U, U \rangle_F^\C - \,{\rm log}\, {\rm r}_E \, {\dim}\, U^2$ is positive definite on ${\rm K}_\C^{\leq w}$.\ps
\noindent Choose now $F={\rm F}_\lambda$ (resp. $F={\rm G}_\lambda$ under ${\rm (GRH)}$), and consider the function $f : \R \rightarrow \R$ defined by $f(t)={\rm cosh}(\frac{t}{2})^{-1}$ (resp. $f=1$). By Lemma \ref{limfl}, the quadratic form of assertion (b) above on ${\rm K}_\C^{\leq w}$ converges to $U \mapsto \langle U, U \rangle_f^\C - \,{\rm log}\, {\rm r}_E \, {\dim}\, U^2$ when $\lambda \rightarrow \infty$. By Proposition \ref{crittnh} (i) and Lemma \ref{matcompl}, this latter form is positive definite on ${\rm K}_\C^{\leq w}$ if, and only if, we have the inequality 
$${\rm log}\, 2\pi\, - {\rm log}\, {\rm r}_E > {\rm t}_w^H$$ 
with $H(t)=F(t)e^{-\frac{t}{2}}$. By Definition \ref{deftwrw}, this property is equivalent to ${\rm r}_E < {\rm r}(w)$ (resp. ${\rm r}_E < {\rm r}^\ast(w)$), which holds by assumption on ${\rm r}_E$. As a consequence, the quadratic form $U \mapsto \langle U, U \rangle_{{\rm F}_\lambda}^\C - \,{\rm log}\, {\rm r}_E \, {\dim}\, U^2$ is positive definite on ${\rm K}_\C^{\leq w}$ for $\lambda$ big enough (e.g. by Sylvester's criterion). Fix such a $\lambda$. For each $v \in {\rm S}_\infty(E)$, define ${\rm q}_v$ as $U \mapsto \langle U, U \rangle_{{\rm F}_\lambda}^{E_v} - \,{\rm log}\, {\rm r}_E \, {\dim}\, U^2$. Set also $C \,=\, 2 \,\widehat{{\rm F}_\lambda}(\frac{i}{4\pi})$ (a positive real number by Footnote \ref{footis4p}). The assertion of the statement holds, by the equivalence of (a) and (b) above.  \end{pf}

\noindent We now prove Theorem \ref{mainthm} of the introduction.\ps

\begin{pf} (of Theorem \ref{mainthm}) Fix a number field $E$ and an integer $w\geq 0$ as in the statement of Theorem \ref{mainthm}. Let $\pi$ be a cuspidal algebraic automorphic representation of ${\rm GL}_n$ over $E$ with weights in $\{0,\dots,w\}$. Then $\pi$ is effective of motivic weight $0 \leq w' \leq w$. As we have ${\rm r}(w') \geq {\rm r}(w)$ it is enough to prove the finiteness of the $\pi$ as in the statement whose motivic weight is exactly $w$. \ps
Fix inner products ${\rm q}_v$ on ${\rm K}_{E_v}^{\leq w}$ for $v \in {\rm S}_\infty(E)$, and a constant $C>0$, as in the statement of Proposition \ref{propmainthm}. By Bushnell and Henniart  \cite[Thm. 1]{bh}, we know that for any integer $n\geq1$ and any cuspidal automorphic representation $\pi$ of ${\rm GL}_n$ over $E$ with conductor ${\mathcal N}(\pi)$, we have the inequality 
$${\rm log}\, \,||\,\mathcal{N}(\pi \times \pi^\vee)\,||\, \,\leq \,(2n-1)\, {\rm log}\,\,||\, {\mathcal N}(\pi)\,||. $$
Applying this to a $\pi$ which is algebraic, effective, of motivic weight $w$, and setting $V_v = {\rm L}(\pi_v \otimes |\det|_v^{-w/2})$ for each $v \in {\rm S}_\infty(E)$, we obtain 
\begin{equation}\label{ineqavecbh}
\sum_{v \in {\rm S}_\infty(E)} {\rm q}_{v}(V_v) \,\leq C\,+ (2n-1) \,{\rm log} \,\,||\, \mathcal{N}(\pi )\,||.
\end{equation}
Fix an arbitrary element $v_0$ in ${\rm S}_\infty(E)$; we have $\dim V_v = n = \dim V_{v_0}$ for each $v$ in ${\rm S}_\infty(E)$.
By the trivial finiteness Lemma \ref{lemmefinitude} below applied to the lattice $L =  \oplus_{v \in {\rm S}_\infty(E)} {\rm K}_{E_v}^{\leq w}$, equipped with the orthogonal sum of the real quadratic forms ${\rm q}_v$, and to the affine function defined by $\varphi(\sum_v x_v)=C+(2\dim x_{v_0}-1) \,{\rm log} \,\,||\, \mathcal{N}(\pi )\,||$, there are only finitely many tuples $(V_v)_{v \in {\rm S}_\infty(E)}$ in $L$ satisfying \eqref{ineqavecbh}. In particular, there are only finitely many possibilities for $n$. We conclude the proof by the classical Lemma \ref{HC}. \end{pf} 

\begin{lemme}\label{lemmefinitude} Let $L$ be a lattice in the Euclidean space $V$, and $\varphi : V \rightarrow \R$ an affine function. The set of $v \in L$ with $v \cdot v \leq \varphi(v)$ is finite.
\end{lemme} 
\ps
\begin{pf} Write $\varphi(v)=w \cdot v + \varphi(0)$ for some $w \in V$, and set $|v|^2=v \cdot v$. The inequality $v \cdot v \leq \varphi(v)$ implies $|v-w/2|^2 \leq |w|^2 /4 +\varphi(0)$.
\end{pf}
\ps
\begin{lemme} \label{HC} Fix $E$ a number field, $\mathcal{N}$ an ideal of $\mathcal{O}_E$, $n\geq 1$ an integer, and for each $v \in {\rm S}_\infty(E)$ an irreducible admissible Harish-Chandra module $U_v$ of ${\rm GL}_n(E_v)$. Then there are only finitely many cuspidal automorphic representations $\pi$ of ${\rm GL}_n$ over $E$, whose conductor is $\mathcal{N}$, and with $\pi_v \simeq U_v$ for each $v \in {\rm S}_\infty(E)$.
\end{lemme}
\ps
\begin{pf} This is a well-known consequence of Harish-Chandra's general finiteness theorem \cite[Thm. 1, p. 8]{harishchandra}, and of the Jacquet, Piatetski-Shapiro, and Shalika, theory of newforms \cite[\S 5, Thm. (ii)]{jpsnew}. 
\end{pf}

\section{Complements: effectiveness and optimality}\label{complements}

\subsection{Effectiveness of the proof of Theorem \ref{mainthm}} Let $w \geq 0$ be an integer, $E$ a number field, and $\mathcal{N} \subset \mathcal{O}_E$ a nonzero ideal. Denote by ${\mathcal S}(E,\mathcal{N},w)$ the set of cuspidal algebraic automorphic representations of general linear groups over $E$, with conductor $\mathcal{N}$, and which are effective of motivic weight $w$. 
Assume ${\rm r}_E < {\rm r}(w)$, so ${\mathcal S}(E,\mathcal{N},w)$ is finite by Theorem \ref{mainthm}. In this generality, it seems very difficult to describe this finite set. Our method provides nevertheless some general effective informations on ${\mathcal S}(E,\mathcal{N},w)$, and even leads to an explicit upper bound on its size in some situations, {\it e.g.} for $||\mathcal{N}||=1$ (case ''unramified at all finite places"). We explain what we mean in the paragraphs (1) and (2) below. \ps

We denote by $K$ the subgroup of $\prod_{v \in {\rm S}_\infty(E)} {\rm K}_{E_v}^{\leq w}$ consisting of the elements  $V=(V_v)$ such that $\dim V := \dim V_v$ is independent of $v$.  Let ${\rm q}_F$ be the quadratic form on $K$ defined by the left-hand side of Formula \eqref{expleffective}, with $F$ as in \S \ref{quadgamma}. For $\varpi$ in ${\mathcal S}(E,\mathcal{N},w)$ we denote by ${\rm V}(\varpi)$ the element $({\rm L}(\varpi_v \otimes |\det |_v^{-\frac{w}{2}}))$ of $K$. \ps\ps

(1) {\it We claim first that there is a computable subset $S'$ of $K$ such that for any $\varpi$ in ${\mathcal S}(E,\mathcal{N},w)$ we have ${\rm V}(\varpi) \in S'$}. For $\lambda>0$ denote by $S_\lambda$ the set of effective $V$ in $K$ satisfying
\begin{equation} \label{discussioneffec}{\rm q}_{{\rm F}_\lambda}(V) \leq  \frac{8}{\pi^2} \lambda \, +\,(\dim V - \frac{1}{2})\, {\rm log}\,|| \mathcal{N}||.\end{equation}
(recall $\widehat{{\rm F}_\lambda}(\frac{i}{4\pi}) = \frac{8}{\pi^2} \lambda$ and ${\rm F}_\lambda(0)=1$) We first choose $\lambda_0>0$ such that ${\rm q}_{{\rm F}_{\lambda_0}}$ is positive definite on $K$. For this step it is important to be able to numerically compute the linear forms ${\rm J}_{{\rm F}_\lambda}^{E_v}$, hence the form ${\rm q}_{{\rm F}_\lambda}$, and we refer to \cite[Prop. IX.3.18]{CL} for a solution to this computational problem. That being done, the ${\rm LLL}$-algorithm\footnote{Use {\it e.g.} the algorithm $\texttt{qfminim}$ in $\texttt{PARI}$, as well as Lemma \ref{lemmefinitude}.} allows to determine the finite set $S_{\lambda_0}$. The set $S'=S_{\lambda_0}$ does the trick, although it is often much too large in practice, and it is better to rather consider $S'=\cap_{\lambda>0} S_\lambda$. This subset of $S_{\lambda_0}$ may be determined by checking again the inequality \eqref{discussioneffec} for each $V \in S_{\lambda_0}$ for other parameters $\lambda$. \ps

\begin{remark} \label{compformql}{\rm We have ${\rm g}'(t)=-\pi (1-t)\, {\rm sin }\, \pi t$ on $[0,1]$, so $g$ is nonincreasing on $[0,+\infty[$ and we have ${\rm F}_\lambda \leq {\rm F}_{\lambda'}$ for $0<\lambda \leq \lambda'$. By Proposition \ref{inegqh1h2} this shows that the $1$-parameter family of quadratic forms ${\rm q}_{{\rm F}_\lambda}$ on $K$ is nondecreasing with $\lambda$; it also converges to ${\rm q}_F$ with $F(t)={\rm cosh}(\frac{t}{2})^{-1} e^{-\frac{t}{2}}$ by Lemma \ref{limfl}. As the term $\frac{8}{\pi^2} \lambda$ goes to $+\infty$ with $\lambda$, the inequality \eqref{discussioneffec} tends to be stronger for small $\lambda$. }
\end{remark}

(2) Fix $V$ in $S'$. It remains to give an explicit upper bound for the number ${\rm m}(V)$ of elements $\varpi$ in ${\mathcal S}(E,\mathcal{N},w)$ with ${\rm V}(\varpi)=V$. Following an argument essentially due to O. Ta\"ibi, {\it we claim that there is an effective upper bound on ${\rm m}(V)$ if the following inequality holds:
\begin{equation} \label{condexplbound} {\rm q}_{F}(V) >  \,(\dim V - \frac{1}{2})\, {\rm log}\,|| \mathcal{N}||, \, \, \, {\rm with}\,\,\, F(t)=e^{-\frac{t}{2}}{\rm cosh}(t/2)^{-1}.\end{equation}}\par
\noindent Note that this is a non trivial problem {\it a priori}, as Harish-Chandra's finiteness result in Lemma \ref{HC} is not effective. We first choose $\lambda>0$ big enough so that the strict inequality above still holds with $F$ replaced by ${\rm F}_\lambda$, and set $\delta_\lambda(V) = {\rm q}_{{\rm F}_\lambda}(V) -  \,(\dim V - \frac{1}{2})\, {\rm log}\,|| \mathcal{N}||$. We have $\delta_\lambda(V)>0$ by assumption and we claim that we have the inequality 
\begin{equation} \label{inegtaibi} {\rm m}(V) \leq \frac{8\, \lambda }{\delta_\lambda(V)\, \pi^2}.\end{equation}
To prove this, assume $\varpi_1,\varpi_2,\dots,\varpi_r$ are distinct elements of ${\mathcal S}(E,\mathcal{N},w)$, with ${\rm V}(\varpi_i)=V$ for each $i$. Set $n=\dim V$. We apply Mestre's explicit formula to the Rankin-Selberg ${\rm L}$-function $\Lambda(s,\varpi \times \varpi^\vee)$ of the isobaric automorphic representation $\varpi=\varpi_1 \boxplus \varpi_2 \boxplus \dots \boxplus \varpi_r$ of ${\rm GL}_{nr}$ over $E$. The only difference with the case of a cuspidal $\varpi$ recalled in \S \ref{sectrs} is that the ${\rm L}$-function $\Lambda(s,\varpi \times \varpi^\vee)$ now has a pole of order $r$ at $s=0$ and $s=1$. In particular, the inequality of Proposition \ref{inegfex} still holds for $\varpi$ but with the term $\widehat{F}(\frac{i}{4\pi})$ replaced by $r \widehat{F}(\frac{i}{4\pi})$. We also  have :\ps\ps

\noindent (i) \, ${\rm L}(\varpi_v \otimes |\det|_v^{-\frac{w}{2}})\,=\, r\, V_v$ for each $v \in {\rm S}_\infty(E)$,\ps\ps

\noindent (ii) $\mathcal{N}(\varpi \times \varpi^\vee) = \prod_{1 \leq i,j \leq r} \mathcal{N}(\varpi_i \times \varpi_j^\vee)$, and \cite[Thm. 1]{bh} still shows the inequality ${\rm log} \, ||\mathcal{N}(\varpi_i \times \varpi_j^\vee)|| \leq (2n-1) {\rm log} \, ||\mathcal{N}||$ for all $i,j$. \ps\ps

\noindent This proves the inequality ${\rm q}_{{\rm F}_\lambda}(r\,V) \,\leq \,r\, \widehat{{\rm F}_\lambda}(\frac{i}{4\pi}) \,+\, r^2 \,(n-\frac{1}{2})\,{\rm log}\, ||{\mathcal N}||$. As we have ${\rm q}_{{\rm F}_\lambda}(r\,V)=r^2 \, {\rm q}_{{\rm F}_\lambda}(V)$, dividing that inequality by $r^2$ leads to $\delta_\lambda(V)\, \leq \, \frac{1}{r}\, \widehat{{\rm F}_\lambda}(\frac{i}{4\pi})$ and concludes the proof of \eqref{inegtaibi}. $\square$ \ps\ps\ps

\noindent In the special case $||\mathcal{N}||=1$, the inequality \eqref{condexplbound} holds for all $V$ in $K-\{0\}$ as we assumed ${\rm r}_E < {\rm r}(w)$. We have thus proved the:

\begin{thm}\label{thmeff} Let $E$ be an number field and $w\geq 0$ an integer with ${\rm r}_E < {\rm r}(w)$.  The algorithm described above leads to an effective upper bound for the number of cuspidal algebraic automorphic representations $\pi$ of $\GL_n$ over $E$ {\rm (}with $n$ varying{\rm )} whose weights are in $\{0,\dots,w\}$, and such that $\pi_v$ is unramified for each finite place $v$ of $E$.
\end{thm}

We refer to \cite[\S IX.3]{CL} for an application of this method in the special case $E=\Q$ and $w \leq 22$, including a complete description of ${\mathcal S}(\Q,(1),w)$ in this range (which requires many other ingredients!). \ps

\subsection{On the optimality of the assumption ${\rm r}_E < {\rm r}(w)$}\label{chaptopt}

A key ingredient in the proof of Proposition \ref{propmainthm} (hence in that of Theorem \ref{mainthm} as well), is the fact that for $\lambda$ big enough, and $r < {\rm r}(w)$, the quadratic form $U \mapsto \langle U,U\rangle_{{\rm F}_\lambda}^\C- {\rm log}\, r \,\,(\dim U)^2$ is positive definite on ${\rm K}_\C^{\leq w}$. It is legitimate to ask if the choice of other test functions could allow to improve these theorems. The next proposition shows that it is not the case.  \ps

\begin{prop}\label{goodtestfunctions} Let $w\geq0$ be an integer, $r \geq 1$ a real number, and set $E=\R$ or $E=\C$. There is an equivalence between:\ps
\begin{itemize} 
\item[(i)] there is a nonnegative test function $F$ satisfying {\rm (POS)} and such that the quadratic form $U \mapsto \langle U,U\rangle_F^E- F(0)\frac{[E:\R]}{2}{\rm log}\, r \,\,(\dim U)^2$ is positive definite on ${\rm K}_E^{\leq w}$,\ps
\item[(ii)] the quadratic form $U \mapsto \langle U,U\rangle_{{\rm F}_\lambda}^E- \frac{[E:\R]}{2}{\rm log}\, r \,(\dim U)^2$ is positive definite on ${\rm K}_E^{\leq w}$ for $\lambda$ big enough,\ps
\item[(iii)] the inequality $r< {\rm r}(w)$ holds.
\end{itemize}
More precisely, if we have $r>{\rm r}(w)$, then there is an effective element $U$ in ${\rm K}_E^{\leq w}$ such that for {\bf any} nonnegative test function $F$ satisfying {\rm (POS)} we have $\langle U,U\rangle_{F}^E- F(0) \frac{[E:\R]}{2}\,{\rm log}\, r \,\,(\dim U)^2 < 0$.
\end{prop}
\ps
\noindent This proposition admits the following variant.
\ps

\begin{prop} \label{goodtestfunctionsGRH} The statements of Proposition \ref{goodtestfunctions} still hold if we replace everywhere ``{\rm (POS)}'' by ``$\widehat{F}(\xi) \geq 0$ for all $\xi \in \R$'', and ${\rm r}(w)$ by ${\rm r}^\ast(w)$.
\end{prop}

\begin{pf} We will prove simultaneously Propositions \ref{goodtestfunctions} ({\it case $A$}) and \ref{goodtestfunctionsGRH} ({\it case $B$}). We first show the equivalence between properties (i), (ii) and (iii). Note that (i) (resp. (ii)) holds for $E=\R$ if, and only if, it holds for $E=\C$, by Proposition \ref{wrwc}. As property (iii) does not depend on $E$, we may and do assume that we have $E=\C$. The implication (ii) $\Rightarrow$ (i) is obvious. Moreover, the implication (iii) $\Rightarrow$ (ii) has already been explained during the proof of Proposition \ref{propmainthm} (apply Lemma \ref{limfl}, Proposition \ref{crittnh} (i) and Lemma \ref{matcompl}). We are left to show (i) $\Rightarrow$ (iii).\ps
 Let $F$ be a nonzero nonnegative test function $F$ satisying (POS) in {\it case $A$}, $\widehat{F}(\xi) \geq 0$ for all $\xi \in \R$ in {\it case $B$}.  We have $F(0)>0$ by Lemma \ref{lemmainegF}, so by replacing $F$ by $F/F(0)$ we may assume $F(0)=1$ without loss of generality. Set $H(t)=F(t)e^{-\frac{t}{2}}$ for $t\geq 0$. Assume that the quadratic form $U \mapsto \langle U,U\rangle_F^E-{\rm log}\, r \,(\dim U)^2$ is positive definite on ${\rm K}_\C^{\leq w}$. By Lemma \ref{matcompl}, this assumption on $F$ asserts that the quadratic form ${\rm log} \frac{2\pi}{r} \phi_w^2 - {\rm q}_w^H$ is positive definite on $\R^{w+1}$. On the other hand, in both cases Lemma \ref{lemmainegF} imply that we have $H \leq H_0$ with $H_0(t)={\rm cosh}(\frac{t}{2})^{-1} e^{-\frac{t}{2}}$ in {\it case $A$}, and $H_0(t)=e^{-\frac{t}{2}}$ in {\it case $B$}. By Proposition \ref{inegqh1h2}, we have thus the inequality
\begin{equation} \label{inegopt} {\rm log} \frac{2\pi}{r} \phi_w^2 - {\rm q}_w^H \leq {\rm log} \frac{2\pi}{r} \phi_w^2 - {\rm q}_w^{H_0}, \, \, \, \, {\rm for}\, \, {\rm all}\,\, r>0\,\, {\rm and}\, \, {\rm all}\,\, w \geq 0.\end{equation}
This shows that ${\rm log} \frac{2\pi}{r} \phi_w^2 - {\rm q}_w^{H_0}$ is positive definite. But this is equivalent to assertion (iii) by Definition \ref{deftwrw} and Proposition \ref{crittnh} (i) applied to $H_0$, and we are done. \ps

We now refine this analysis to prove the second part of the proposition. We assume $r>{\rm r}(w)$ in {\it case $A$}, and  $r > {\rm r}^\ast(w)$ in {\it case $B$}. Define $H_0$ accordingly as above. We will show that there is an explicit nonzero element $x=(x_k)_{0 \leq k \leq w} \in \R^{w+1}$ satisfying the following two properties:\ps

\noindent (a) $x_k$ is a nonnegative integer, with $x_k = x_{w-k}$ for all $0 \leq k \leq w$,\ps
\noindent (b) the quadratic form ${\rm log} \frac{2\pi}{r} \phi_w^2 - {\rm q}_w^{H_0}$ is $<0$ on $x$.\ps

Note that for any test function $F$ satisfying {\rm (POS)} in {\it case $A$}, or $\widehat{F}(\xi)\geq 0$ for all $\xi \in \R$ in {\it case $B$}, and if we set $H(t)=F(t)e^{-\frac{t}{2}}$, the quadratic form ${\rm log} \frac{2\pi}{r} \phi_w^2 - {\rm q}_w^H$ is negative as well on $x$, by property (b) and the inequality \eqref{inegopt}. Moreover, by the first part of property (a) the element $U= \sum_{k=0}^w x_k \eta^{-w+2k}$ of ${\rm K}_\C^{\leq w}$ has the required properties in the case $E=\C$ by Lemma \ref{matcompl}. By the second part of property (a), $2U$ has the form ${\rm res} \,U'$ for a unique $U' \in {\rm K}_\R^{\leq w}$ with $U' = U' \cdot \epsilon_{\C/\R}$, and $U'$ has the required properties in the case $E=\R$ by the claim of the third paragraph of the proof of Proposition \ref{wrwc}.\ps

We are thus left to prove the existence of nonzero $x \in \R^{w+1}$ satisfying (a) and (b) above. By the assumption $r>{\rm r}(w)$ in {\it case $A$}, and  $r > {\rm r}^\ast(w)$ in {\it case $B$}, the nonzero element $y={\rm v}_w^{H_0}$ satisfies (b), as well as $y_{w-k}=y_k$ for all $0 \leq k \leq w$ by Proposition \ref{orthpol} (this is equivalent to say that ${\rm P}_w^H$ is real). So a suitable integer multiple of $y$ does the trick if we know that its coefficients are nonnegative rational numbers. This concludes the proof in {\it case $B$}, as we have $4^w \,y_k\,=\, {{2k}\choose{k}} {{2(w-k)}\choose{w-k}}$ by Proposition \ref{calctnGRH} and Formula \eqref{explpnh0}. \ps

Suppose now that we are in {\it case $A$}. As explained in the paragraph preceding Corollary \ref{calctnSANSGRH}, we expect that ${\rm v}_w^{H_0}$ has positive coefficients for all $w$; we know it by explicit computations for small $w$, and in particular for $w \leq 25$. The element ${\rm v}_w^{H_0}$ has rational coefficients for $w$ odd by Corollary \ref{calctnSANSGRH}, but not for $w$ even. This is not a problem, since assuming ${\rm v}_w^{H_0}$ has nonnegative coefficients we may always choose an element $y \in \Q^{w+1}$, with $y_k=y_{w-k}$ for all $0 \leq k \leq w$, with nonnegative coefficients, which is close enough to ${\rm v}_w^{H_0}$ so that the quadratic form ${\rm log} \frac{2\pi}{r} \phi_w^2 - {\rm q}_w^{H_0}$ is $<0$ on $y$. This concludes the proof of the proposition for $w \leq 25$. The remaining case $w > 25$ actually follows from this one. Indeed, we have then ${\rm r}(w) \leq {\rm r}(25) < 1 \leq r$ by Corollary \ref{corthsansgrh}. We conclude as ${\rm K}_\C^{\leq w}$ contains ${\rm K}_\C^{\leq 25}$ for $w$ odd, ${\rm K}_\C^{\leq 24}$ for $w$ even.\end{pf}

\begin{example} \label{exampleQ} {\rm Let us discuss the optimality of Theorem \ref{thmch}. Take $H(t)=e^{-t/2}$ and consider the element $x=(4,1,1,\dots,1,1,4) \in \R^{26}$, which generates a line close to $\R {\rm v}_{25}^H$. The quadratic form ${\rm log}\, 2\pi \,\,\phi_{25}^2 - {\rm q}_{25}^H$ is negative at $x$ (it is $\approx -1,04$ up to $10^{-2}$).  Set $$U={\rm I}_1+{\rm I}_3 + \cdots + {\rm I}_{21}+{\rm I}_{23} + 4 \,{\rm I}_{25} \in {\rm K}_\R^{\leq 25}.$$
\noindent By the proof above, we have $\langle U, U \rangle^\R_F <0$ for any nonzero test function $F \geq 0$ satisfying $\widehat{F} \geq 0$. As a consequence, even under {\rm (GRH)}, our method cannot rule out the possibility that there exist infinitely many cuspidal algebraic automorphic representations $\varpi$ of ${\rm GL}_n$ over $\Q$ (with $n$ varying), which are effective of motivic weight $25$, of same Artin conductor, and with ${\rm L}(\varpi_\infty \otimes |\det|^{-\frac{25}{2}}) \in \Z U$. Moreover, if we set 
$$V={\rm I}_2+{\rm I}_4 + \cdots + {\rm I}_{20}+{\rm I}_{22} + 4\, {\rm I}_{24}\in {\rm K}_\R^{\leq 24},$$
Corollary \ref{corthsansgrh} and the proof above show that we have $\langle V, V \rangle^\R_F <0$ for any nonzero test function $F \geq 0$ satisfying (POS).  As a consequence, without {\rm (GRH)}, our method cannot rule out the possibility that there exist infinitely many cuspidal algebraic automorphic representations $\varpi$ of ${\rm GL}_n$ over $\Q$ (with $n$ varying), which are effective of motivic weight $24$, of same Artin conductor, and with ${\rm L}(\varpi_\infty \otimes |\det|^{-12}) \in \Z V$.}
\end{example}

{\scriptsize 
\subsection{Comparison with the ``standard ${\rm L}$-function'' method}\label{retoursurgj} 

As mentionned in the introduction, applying the explicit formula to Godement-Jacquet ${\rm L}$-functions, rather than to Rankin-Selberg ${\rm L}$-functions, also proves some variant of Theorem \ref{mainthm}. We only briefly discuss it here, and use a small font, because the resulting statements are much weaker. Let $\varpi$ be a cuspidal algebraic automorphic representation of  ${\rm GL}_n$ over $E$ which is effective of motivic weight $w\geq 0$, with conductor $\mathcal{N}$, and with $\varpi \neq 1$. If we denote by $V_v$ the element ${\rm L}(\varpi_v \otimes |.|_v^{-w/2})$ of ${\rm K}_{E_v}^{\leq w}$ for $v \in {\rm S}_\infty(E)$, the basic inequality takes the form
{\scriptsize \begin{equation}\label{gjineg} \sum_{v \in {\rm S}_\infty(E)} \left({\rm J}_{F}^{{\rm E}_v}(V_v) \,-\, F(0) \,\frac{[E_v:\R]}{2}\,{\rm log}\, {\rm r}_E\,\dim V_v \right) \, \leq \,\frac{F(0)}{2}\, {\rm log} \,||\mathcal{N}||.\end{equation}}
\par \noindent Here we have to assume that $F$ satisfies {\rm (POS)} and that its support is in $[-{\rm log}\, 2, {\rm log}\, 2]$ since the analogue of \eqref{posanpipiprime} does not hold anymore. To fix ideas we just choose $F={\rm F}_{{\rm log}\, 2}$.  
In order to draw any conclusion on $||\mathcal{N}||$ the sum on the left-hand side of \eqref{gjineg} has to be $>0$. The sequence of real numbers $({\rm J}_{{\rm F}_{{\rm log}\, 2}}^\C(\eta^w))_{w \geq 0}$ is decreasing, and a numerical computation shows that for $w\geq 2$ it is smaller than $2 {\rm J}_F^\R(\epsilon_{\C/\R})$. This shows that for $E_v=\C$, the $v$-part of the sum in \eqref{gjineg} is bounded below by $n({\rm J}_F^\C(\eta^w) -\,{\rm log}\, {\rm r}_E)$, and that for $E_v=\R$ this $v$-part is bounded below by $\frac{n}{2}({\rm J}_F^\C(\eta^w) -\,{\rm log}\, {\rm r}_E)$ for $w>0$, and by $\frac{n}{2}(2 {\rm J}_F^\R(\epsilon_{\C/\R}) -\,{\rm log}\, {\rm r}_E)$ for $w=0$. \ps

For $w \geq 0$ an integer we set ${\rm s}(w)={\rm exp}\, ({\rm J}_{{\rm F}_{{\rm log}\, 2}}^\C(\eta^w))$. This is a decreasing sequence, and we have ${\rm s}(10) \approx 1.003$ up to $10^{-3}$ and ${\rm s}(w)<1$ for $w>10$. Table \ref{tablesw} gives the relevant numerical values of ${\rm s}(w)$. We also set ${\rm s}'(0)={\rm exp}\,( 2\, {\rm J}_{{\rm F}_{{\rm log}\, 2}}^\R(\epsilon_{\C/\R}))$, a real number which is $\approx 2.323$ up to $10^{-3}$. The analysis above proves a variant of Theorem \ref{mainthm} in which ${\rm r}(w)$ is replaced with ${\rm s}(w)$, unless we have $w=0$ and $E$ has a real place in which case ${\rm r}(0)$ has to be replaced with ${\rm s}'(0)$. They are much lower than ${\rm r}(w)$. Actually, there are only five number fields $E$ with ${\rm r}_E <  {\rm s}(0) \approx 2.669$, namely $\Q$ and $\Q(\sqrt{d})$ for $d=-3, -4, 5$ and $-7$. Let us add that we would not gain anything using {\rm GRH} here.

\begin{table}[!h]
\renewcommand{\arraystretch}{1.5}
{\scriptsize
{
\begin{tabular}{c||c|c|c|c|c|c|c|c|c|c|c}
 $w$  & $0$ & $1$ & $2$ & $3$ & $4$ & $5$ & $6$ & $7$ & $8$ & $9$ & $10$ \\ 
\hline ${\rm s}(w)$ & $2.67$ & $2.34$ & $2.06$ & $1.84$ & $1.66$ & $1.50$ & $1.37$ & $1.26$ & $1.16$ & $1.08$ & $1.00$ \\
\end{tabular}\ps\ps
\caption{{\small Values of ${\rm s}(w)$ up to $10^{-2}$.}}
\label{tablesw}
}
}\end{table}
}

\end{document}